\documentclass[12pt]{amsart}
\usepackage{amsmath , amsthm , amsfonts , ifpdf}
\usepackage{color}
\usepackage{graphicx}
\usepackage{multicol}
\usepackage{epstopdf}
\usepackage{verbatim}

\usepackage{tikz}

\usepackage{amsmath}
\usepackage{amssymb}

\usetikzlibrary{arrows , shapes , trees , backgrounds}
\usetikzlibrary{intersections}

\newif\ifdebug 

\newcommand{\norm}[1]{{\left | \left |  #1  \right | \right |}}
\newcommand{\inner}[2]{{\left < #1, #2 \right >}}

\renewcommand{\div}{{{div}}}
\newcommand{\spt}{{\text{spt}}}

\newcommand{\diam}{\text{diam}}
\newcommand{\Tr}{\mathcal{T}}

\theoremstyle{plain}
\newtheorem{theorem}{Theorem}[section]
\newtheorem*{theorem*}{Theorem}

\newtheorem{pro}[theorem]{Proposition}
\newtheorem{Def}[theorem]{Definition}
\newtheorem{lem}[theorem]{Lemma}

\newtheorem{cor}[theorem]{Corollary}

\theoremstyle{definition}
\newtheorem*{Def*}{Definition}

\newtheorem{Rem}[theorem]{Remark}

\numberwithin{equation}{section}

\newcommand{\bpo}{\begin{pro}}
	\newcommand{\epo}{\end{pro}}
\newcommand{\be}{\begin{equation}}
\newcommand{\ene}{\end{equation}}
\newcommand{\br}{\begin{Rem}}
	\newcommand{\er}{\end{Rem}}
\newcommand{\bl}{\begin{lem}}
	\newcommand{\el}{\end{lem}}
\newcommand{\bd}{\begin{Def}}
	\newcommand{\ed}{\end{Def}}
\newcommand{\ben}{\begin{enumerate}}
	\newcommand{\een}{\end{enumerate}}
\newcommand{\bp}{\begin{proof}}
	\newcommand{\ep}{\end{proof}}
\newcommand{\beq}{\begin{equation*}}
\newcommand{\eeq}{\end{equation*}}
\newcommand{\bear}{\begin{eqnarray*}}
	\newcommand{\eear}{\end{eqnarray*}}
\newcommand{\bt}{\begin{theorem}}
	\newcommand{\et}{\end{theorem}}
\newcommand{\bst}{\begin{split}}
	\newcommand{\est}{\end{split}}

\newcommand{\bal}{\begin{aligned}}
	\newcommand{\eal}{\end{aligned}}

\renewcommand{\P}{\partial}
\newcommand{\F}[2]{\frac{#1}{#2}}
\newcommand{\la}{\langle}
\newcommand{\ra}{\rangle}

\newcommand{\PLH}{{\mkern-1mu\times\mkern-1mu}}
\newcommand{\R}{\mathbb{R}}
\newcommand{\AF}{\mathfrak{F}_{\phi}}
\newcommand{\bM}{\mathbb{M}}
\newcommand{\mG}{\mathcal{G}}
\newcommand{\IMC}{integer multiplicity current}

\newcommand{\Ca}{Caccioppoli}
\newcommand{\Sc}{\varepsilon}
\newcommand{\wrt}{with respect to}
\newcommand{\vp}{\varphi}
\newcommand{\nb}{\nabla}

\begin{document}
	
	\title[Curvature Estimates]{The area minimizing problem in conformal cones: II}
	\author{Qiang Gao, Hengyu Zhou}
	\date{\today}
\address[Q.~G]{Department of Mathematics, Sun Yat-Sen University, 510275, Guangzhou, P. R. China}
\email{gaoqiangks@outlook.com}

\address[H. ~Z]{ College of Mathematics and Statistics, Chongqing University, Huxi Campus, Chongqing, 401331, P. R. China}
\address{Chongqing Key Laboratory of Analytic Mathematics and Applications, Chongqing University, Huxi Campus, Chongqing, 401331, P. R. China}
\email{zhouhyu@cqu.edu.cn}
	\date{\today}
\begin{abstract}
	In this paper we continue to study the connection among  the area minimizing problem, certain area functional and the Dirichlet problem of minimal surface equations in a class of conformal cones with a similar motivation from \cite{GZ20}. These cones are certain generalizations of hyperbolic spaces. We describe the structure of area minimizing $n$-{\IMC}s in bounded $C^2$ conformal cones with prescribed $C^1$ graphical boundary via a minimizing problem of these area functionals. As an application we solve the corresponding Dirichlet problem of minimal surface equations under a mean convex type assumption. We also extend the existence and uniqueness of a local area minimizing {\IMC} with star-shaped infinity boundary in hyperbolic spaces into a large class of complete conformal manifolds. 
\end{abstract}
\maketitle
\section{Introduction} 
 In this paper we continue to study the area minimizing problem with prescribed boundary in a class of conformal cones similar to \cite{GZ20}.  A conformal cone in this paper is defined as follows. 
\bd\label{def:Model}
Let $N$ be a $n$-dimensional Riemannian manifold with a metric $\sigma$, let $\R$ be the real line with the metric $dr^2$ and $\phi(x)$ be a $C^3$ positive function on $N$. In this paper we call  
\be \label{eq:Model}
M_{\phi}:=\{N\PLH\R,  \phi^2(x)(\sigma+dr^2) \}
\ene
as a conformal product manifold. Let $\Omega$ be a $C^2$ bounded domain in $N$. We refer $\Omega\PLH \R$ in $M_\phi$ as a conformal cone, denoted by $Q_\phi$. 
\ed 
 Let $\psi (x)$ be a $C^1$ function on $\P\Omega$ and $\Gamma$ be its graph in $\P\Omega\PLH\R$.  The area minimizing problem in a conformal cone $Q_\phi$ is to find a $n$-{\IMC} in $\bar{Q}_\phi$, the closure of $Q_\phi$, to realize 
\be \label{mass:minimizing:A}
\min\{\bM(T) \mid T\in  \mG, \P T=\Gamma \}
\ene
where $\bM$ is the mass of {\IMC}s in $M_\phi$, and $\mG$ denotes the set of n-{\IMC}s with compact support in $\bar{Q}_\phi$, i.e. for any $T\in \mG$, its support $spt(T)$ is contained in $\bar{\Omega}\PLH [a,b]$ for some finite numbers $a<b$. See subsection \ref{sec:IMC} for more details. \\
\indent The main reason to study the conformal product manifold $M_\phi$ in Definition \ref{def:Model} is that hyperbolic space is a special case of $M_\phi$ (see remark \ref{Rm:hyperbolic:space}). With this model in \cite[Theorem 4.1]{HL87}, Hardt-Lin showed that there is a unique local area minimizing $n$-{\IMC} for any prescribed infinity $C^1$ star-shaped boundary. Moreover it is a radial minimal graph over $S^n_+$ in hyperbolic spaces. On the other hand in \cite[Section 4.1]{Lin85}, Lin described the solution to the area minimizing problem \eqref{mass:minimizing:A} with $C^1$ graphical boundary in a bounded $C^2$ cylinder via a minimizing problem of an area functional of bounded variation (BV) functions. Motivated by Lin's idea in \cite{GZ20} we studied the area minimizing problem \eqref{mass:minimizing:A} in a conformal product manifold $M_h =\{N\PLH\R, h^2(r)(\sigma+dr^2)\}$. Note that there $h(r)$ only depends on $r\in \R$. Based on these preceding results, it is natural to consider the corresponding area minimizing problem \eqref{mass:minimizing:A} in conformal cones of $M_\phi$. In particular it is desirable to know how many phenomena in hyperbolic sapce does not really rely on the hyperbolic structures. We refer readers to \cite{Alm87} for some history remarks on such kinds of area minimizing problems. \\
\indent  For any open set $W$ we denote the set of all bounded variation functions on $W$ by $BV(W)$.  A key concept for our study of the problem \eqref{mass:minimizing:A} is an area functional in  $BV(W)$ defined as follows. 
\begin{equation} \label{def:functional:A}
	\begin{split}
		\AF(u, W) &= \sup \{\int_\Omega \phi^n(x)h + u \div(\phi^n(x)X) dvol \mid h \in C_0(W),\\
					   &X \in T_0(W) \text{ and } h^2+\inner{X}{X} \le 1\}
	\end{split}
\end{equation}
where $dvol$ and $div$ are the volume form and the divergence of $N$ respectively, $C_0(W)$ and $T_0(W)$ denote the set of smooth functions and vector fields with compact support in $W$ respectively. Note that  when $u \in C^1(W)$, $\AF(u,W)$ is the area of the graph of $u(x)$ in $M_\phi$. \\
\indent  Let $\Omega$ be the $C^2$ domain in Definition \ref{def:Model}. Suppose $ \Omega'$ is a $C^2$ domain in $N$ such that $\Omega\subset\subset \Omega'$. That is, the closure of $\Omega$ is a compact set in $\Omega'$. Suppose $\psi(x)\in C^1(\Omega'\backslash \Omega)$. The following minimizing problem 
\be \label{BV:minimizing:functional}
\min\{ \AF(v, \Omega') \mid v(x)\in BV(\Omega'), v(x)=\psi(x) \text{  on } \Omega'\backslash \Omega\}
\ene 
also plays an important role to solve \eqref{mass:minimizing:A}. \\
\indent In fact the key idea to solve \eqref{mass:minimizing:A} is to establish the connection among the problem \eqref{mass:minimizing:A}, the area functional minimizing problem \eqref{BV:minimizing:functional} and the Dirichlet problem of minimal surface equations in $M_\phi$. This can be easily seen when $\Sigma$ is a minimal graph of $u(x)$ in $M_\phi$ over $\Omega$ with $C^1$ boundary $\psi(x)$ on $\P\Omega$. From Theorem \ref{thm:minimal:equation} one has (1) $u(x)$ should satisfy 
\be \label{eq:minimal:equation:Mphi}
	div(\F{Du}{\omega})+n\la D\log\phi, \F{Du}{\omega}\ra=0  \text{ on } \Omega
\ene 
with $u(x)=\psi(x)$ on $\P \Omega$;
 (2) the area of $\Sigma$ is less than that of any compact $C^2$ surface $S$ with $\P  S$ as the graph of $\psi(x)$ containing in $\bar{\Omega}\PLH\R$; (3) $u(x)$ realizes the minimum of \eqref{BV:minimizing:functional} if requiring all $v(x)\in C^2(\Omega)\cap BV(\Omega')\cap C(\bar{\Omega})$. \\
\indent In $M_\phi$ a similar connection between the problem \eqref{mass:minimizing:A} and the problem \eqref{BV:minimizing:functional} is obtained in Theorem \ref{thm:area:minimizing}. It says that if $u(x)$ is the solution to the problem \eqref{BV:minimizing:functional}, then $T=\P [[U]]|_{\bar{Q}_\phi}$ solves the problem \eqref{mass:minimizing:A} in $M_\phi$ where $U$ is the subgraph of $u(x)$ and $[[U]]$ is the corresponding {\IMC}. This generalizes Lin's result \cite{Lin85} into conformal cones in $M_\phi$. Its proof is based on the following three observations: (1) by Theorem \ref{subgraph_and_area}, $\AF(u,.)$ is just the perimeter of its subgraph $U$; (2) if $T$ solves \eqref{mass:minimizing:A} and $\Gamma=(x,\psi(x))$ for $\psi(x)\in C^1(\P\Omega)$, then $T$ is a boundary of a {\Ca} set (see Lemma \ref{au:mass:minimizing}); (3) by Theorem \ref{thm:decrease_perimeter}, for any {\Ca} set $F$ in $Q\subset M_\phi$ with $\P F\subset \Omega\PLH [a,b]$, there is a $u(x)\in BV(\Omega)$ such that $\AF(u,\Omega)$ is less than the perimeter of $F$. Note that the observation (3) is a very general phenomenon in area-type functionals  \cite{Mir64}, \cite{Lin89}, \cite{NS96}, \cite{Giu84},\cite{Ger98} and \cite{BD87} (see remark \ref{Rm:decreasing:perimeter}). We refer it as the Miranda's observation.  \\
\indent  A direct application of Theorem \ref{thm:area:minimizing}  is the Dirichlet problem of minimal surface equations in $M_\phi$.  In Theorem \ref{thm:DR:phimeanconvex} we show that if $\Omega$ is $\phi$-mean convex, then the Dirichlet problem \eqref{eq:minimal:equation:Mphi} with continuous boundary data has a unique solution in $C^2(\Omega)\cap C(\bar{\Omega})$. Here $\Omega$ is $\phi$-mean convex if the mean curvature of $\P\Omega$ satisfies that $
H_{\P\Omega}+n \la \vec{\gamma}, D\log\phi \ra\geq 0 \text{ on } \P\Omega$
where $\vec{\gamma}$ is the outward normal vector of $\P\Omega$ and $H_{\P\Omega}=div(\vec{\gamma})$. A more general form of the Dirichlet problem \eqref{eq:minimal:equation:Mphi} was already considered by Casteras-Heinonen-Holopainen \cite{CHH19}.  To obtain the $C^0$ esimate they relied on a lower bound of the Ricci curvature $\Omega$ and the positive mean curvature of $\P\Omega$ (see remark \ref{Rm:CHHwork}) . But out result in Theorem \ref{thm:DR:phimeanconvex} is independent of the curvature of $\Omega$ and also can not be obtained by the classical continuous method in  \cite[Chapter 11, 17.2]{GT01} (see remark \ref{rk:continuous:method}). A consequence of Theorem \ref{thm:DR:phimeanconvex} is that if $\Omega$ is $\phi$-mean convex, the area minimizing {\IMC} in \eqref{mass:minimizing:A} is unique as the graph of a $C^2$ function to solve  the Dirichlet problem (see Theorem \ref{thm:uniqueness:mass:minimizng}).  \\
\indent At last we consider the existence and unqiueness of local area minimizing {\IMC} in $M_\phi$ with infinity boundary $\Gamma$ when $\phi(x)$ can be written as $\phi(d(x,\P N))$ which  goes to $+\infty$ as $d(x,\P  N)$ goes to zero in $N$. Here $N$ is a compact Riemannian manifold with $C^2$ boundary  and $d$ is the distance function in $N$. Let $N_r$ be the set $\{x\in N: d(x,\P N)>r\}$. In Theorem \ref{thm:generalization:Hyperbolic} we obtain that if there is $r_1$ such that for any $r\in (0, r_1)$ $N_r$ is $\phi$-mean convex, then for any $\psi(x)\in C(\P N)$ and $\Gamma =(x,\psi(x))$ there is a unique local area minimizing {\IMC} $T$ with infinity boundary $\Gamma$. Moreover $T$ is a minimal graph in $M_\phi$ over $N$. This illustrates that Theorem 4.1 in \cite{HL87} does not depend on the hyperbolic structure (see remark \ref{Rm:hyperbolic:space} ). Unlike the case in hyperbolic spaces our existence in Theorem \ref{thm:generalization:Hyperbolic} is from the Dirichlet problem of minimal surface equations  in Theorem \ref{thm:DR:phimeanconvex}, not from the results in \cite{And82,And83} by Anderson via geometric measure theory.   \\
\indent  For the idea of radial graphs in hyperbolic space  we refer to \cite{GS00, DS09,LX12,NS96,Spr07, RT19,GS05,GS08, Lin89} etc. For the Dirichlet problem of minimal surface equations in Riemannian manifolds, we refer to \cite{CHH19}, \cite{AS19}, \cite{Sim97}, \cite{ARS16}, \cite{AY19} and references therein. As for the variational method to study the Dirichlet problem of minimal surface equations, we refer to \cite{Giu84,Giu80,BD87}, \cite{Zhou19} and references therein.  \\
\indent Our paper is organized as follows. In section 2 we show three properties of $C^2$ minimal graphs in $M_\phi$ and collect preliminary facts on BV functions. In section 3 we show the $C^\infty$ approximation theorem of $\AF(u,\Omega)$ (Theorem \ref{approximation_of_area}) and the Miranda's observation (Theorem \ref{subgraph_and_area} and Theorem \ref{thm:decrease_perimeter}). In section 4 we consider the connection between the problem \eqref{mass:minimizing:A} and the problem \eqref{BV:minimizing:functional} (Theorem \ref{thm:area:minimizing}). In section 5 we discuss the Dirichlet problem of minimal surface equations in $M_{\phi}$ on a $\phi$-mean convex domain (Theorem \ref{thm:DR:phimeanconvex}). As an application we obtain the uniqueness of the problem \eqref{mass:minimizing:A} under the $\phi$-mean convex assumption (Theorem \ref{thm:uniqueness:mass:minimizng}). In section 6 we discuss the existence and uniqueness of local area minimizing $n$-{\IMC} with infinity graphical boundary (Theorem \ref{thm:generalization:Hyperbolic}). In Appendix A we give a proof of  the interior estimate of minimal surface equations in \eqref{eq:minimal:equation:Mphi} (Theorem \ref{thm:interior:estimate}) following from an idea of Eichmair \cite{Eich09}. \\
\indent The first author is supported by the National Natural Science Foundation of China, No. 11771456. The second author is supported by the National Natural Science Foundation of China, No.
11801046 and by the Fundamental Research Funds for the Central Universities, China, No.2019CDXYST0015. 
\section{Preliminaries}
Throughout this section we assume $N$ is a manifold with a metric $\sigma$, $\phi(x)>0$ is a positive $C^3$ function defined on $N$ and $M_\phi$ is the conformal product manifold $N\PLH \R $ with the metric $\phi^2(x)(g + dr^2)$. We show three properties of $C^2$ minimal graphs in $M_\phi$ and collect some results on BV functions in Riemannian manifolds for later use. 
\subsection{Three properties}  
\bd Let $S$ be a $C^2$ orientable hypersurface in a Riemannian manifold $M$ with a normal vector $\vec{v}$. We call $div(\vec{v})$ as the mean curvature of $S$ with respect to $\vec{v}$. 
\ed 
 Let $\Omega$ be a bounded domain in $N$ and $u\in C^2(\Omega)$.  Let $\Sigma$ be the graph of $u(x)$. Then its upward normal vector in the product manifold $N\PLH\R$ is $\vec{v}_\Sigma:=\F{\P_r-Du}{\omega}$ where $D$ denotes the gradient of $N$ and $\omega=\sqrt{1+|Du|^2}$. Then we have 
\bl \label{lm:mean:curvature} The mean curvature of $\Sigma$ with respect to its upward normal vector in $M_\phi$ is 
\be 
    H_\Sigma=\F{1}{\phi(x)}(-div(\F{Du}{\omega})-n \la D\log\phi, \F{Du}{\omega}\ra )
\ene 
where $\omega=\sqrt{1+|Du|^2}$, $n$ is the dimension of $N$ and $div$ is the divergence of $N$. 
\el 
\bp  We write the metric of $M_\phi$ as $e^{2f}(\sigma+dr^2)$ where $f=\log\phi$. The upward normal vector of $\Sigma$ in $M_{\phi}$ is $e^{-f} \vec{v}_\Sigma$. 
By  \cite[Lemma 3.1]{Zhou19-3},  the mean curvature of $\Sigma$ in $M_\phi$ with respect to $e^{-f}\vec{v}_\Sigma$ is 
\be \label{det:expression}
H_\Sigma = e^{-f}( H+ndf(\vec{v}_\Sigma))
\ene 
where $H$ is the mean curvature of $S$ with respect to $\vec{v}_\Sigma$ in the product manifold $N\PLH\R$ and $n$ is the dimension of $N$. With a straightforward computation $H$ is $-div(\F{Du}{\omega})$ where $div$ is the divergence of $N$.  Putting $\vec{v}_\Sigma$ and $f=\log\phi$ into \eqref{det:expression}, we obtain the conclusion. 
\ep 
Now we summarize three properties of minimal graphs in $M_\phi$ as follows. 
\bt \label{thm:minimal:equation} Suppose $\Sigma$ is a minimal graph of $u(x)$ over $\Omega$ with $C^1$ boundary $\Gamma:=\{(x,\psi(x)):x\in \P\Omega, \psi(x)\in C^1(\P\Omega)\}$. Then 
\begin{enumerate}
	\item $u(x)$ satisfies that 
\be  \label{eq:minimal:surface:equation:phi}
-div(\F{Du}{\omega})-n \la D\log\phi (x), \F{Du}{\omega}\ra=0
\ene 
\item Let $\Omega'$ be a domain such that $\Omega\subset\subset \Omega'$ and $\psi(x)\in C^1(\Omega'\backslash\Omega)$. Further assuming $u(x)= \psi(x)$ on $\Omega'\backslash \Omega$, then $u(x)$ realizes  
\be 
\min\{ \AF(v, \Omega') \mid v(x)\in C^2(\Omega), v(x)=\psi(x) \text{  on } \bar{\Omega}'\backslash \Omega\}
\ene 
where  $\AF(v,\Omega')=\int_{\Omega'} \phi^{n}(x)\sqrt{1+|Dv|^2}dvol$ and $dvol$ is the volume form of $N$.
\item The area of $\Sigma$ achieves the minimum of the area of all $C^2$ compact orientable hypersurface $S$ in $\bar{Q}_\phi:=\bar{\Omega}\PLH \R$ with $\P S=\Gamma$.
\end{enumerate}
\et 
\bp Note that the  property (1) is from Lemma \ref{lm:mean:curvature}.  \\
\indent For any $C^2$ graph $S=(x,v(x))$ over $\Omega$, its area in $M_\phi$ is 
\begin{equation} \label{def:area}
	\int_{\Omega} \phi^n(x)\sqrt{1+|Dv|^2}dvol
\end{equation}
where $dvol$ is the volume form on $\Omega$. Thus the property (2) follows from the property (3). It is sufficient to show the property (3). \\
\indent Recall that $\vec{v}_\Sigma:=\F{\P_r-Du}{\omega}$. For any $t\in\R$, define a map $T_t:N\PLH\R \rightarrow N\PLH\R$ as  $T_t(x,r)=(x,r-t)$.  We define a new vector field in $\bar{Q}_\phi$ as 
\be 
 X(x, u(x)-t)= \F{1}{\phi(x)}T_{t*}(\vec{v}_\Sigma) 
\ene 
where $T_{t*}$ is the pushforword of $T_t$. 
Thus $X$ is a smooth unit vector field in the tangent bundle of $\bar{Q}_\phi$.  By the definition of divergence, (see Page 423 in \cite{Lee2013}) we define an $n$-form as 
\be 
	\xi = X\lrcorner (\phi^{n+1} dvol\wedge dr) = X \lrcorner dvol_{M_\phi}
\ene

A key fact is 
\be 
d \xi = div_{M_\phi}(X)dvol_{M_\phi}= H_{T_t(\Sigma)}dvol_{M_\phi}
\ene  where $div_{M_\phi}$ and $dvol_{M_\phi}$ is the divergence and volume form of $M_\phi$.  Note that $T_t$ is an isometry of $M_\phi$ for every $t\in\R$. Because $\Sigma$ is minimal in $M_\phi$, so is $T_t(\Sigma)$. Then $d\xi=0$ on the whole $Q_\phi$.  \\
\indent  Now let $S$ be any compact $C^2$ orientable hypersurface in $\bar{Q}_\phi$ satisfying $\P S=\Gamma$. Without loss of generality we can assume $S$ and $\Sigma$ bounds a domain $W$ in $M_\phi$. Suppose that the outward normal vector of $W$ on $\Sigma$ and $S$ are denoted by $\vec{n}_\Sigma$ and $\vec{n}_S$, the volume form of $\partial W$ on $\Sigma$ and $S$ are denoted by $dvol_\Sigma$ and $dvol_S$ respectively.  It's clear that the upward normal vector of $\Sigma$ in $M_\phi$,  $\vec{n}_\Sigma$, is $\frac{1}{\phi}\vec{v}_\Sigma$.  Applying the divergence theorem in the domain $W$, we obtain 
\begin{align*}
0&=\int_W d\xi =\int_\Sigma \inner{X}{\vec{n}_\Sigma}_{M_\phi} dvol_\Sigma - \int_S \inner{X}{\vec{n}_S}_{M_\phi} dvol_S\\
	 & \geq Area(\Sigma)-Area (S)
\end{align*} 
with equality if and only if $S=\Sigma$. 
Thus we obtain the property (3). The proof is complete. 
\ep
\subsection{BV functions} Now we collect some preliminary facts for BV functions in Riemannian manifolds.  We refer readers to the books of Giusti\cite{Giu84}, Evans-Gariepy\cite{EG15} and the papers of  McGonagle-Xiao \cite{MX13}, \cite{Zhou19} etc. Recall that $\Omega$ is a bounded domain in $N$ and $dvol$ is the volume form of $N$.  Let $T_0\Omega$ be the set of smooth vector fields with compact support in $\Omega$.\\
\indent  Now we define BV functions and {\Ca} sets as follows:
\bd
	Let $u \in L^1(\Omega)$.  Define
	\begin{equation}\label{def:bv}
		\norm{Du}(\Omega)=\sup\{\int_\Omega u \div(X)dvol \mid  X\in T_0\Omega \text{ and } \inner{X}{X} \le 1\}
	\end{equation}
	where $dvol$ is the volume form of $\Omega$. 
	If $\norm{Du}(\Omega)<\infty$, we say that $u$ has bounded variation in $\Omega$. The set of all functions with bounded variation is denoted by $BV(\Omega)$. If $u$ belongs to $BV(W)$ for any bounded domain $W\subset \subset \Omega$, we say $u\in BV_{loc}(\Omega)$. \\
\indent Let $E$ be a Borel set in $\Omega$ and $\lambda_E$ be the characteristic function of $E$. If $\lambda_E\in BV_{loc}(\Omega)$,  then $E$ is called as a {\Ca} set and $||D\lambda_E||(\Omega)$ is called as the perimeter of $E$. In the reminder of this paper, we also write it as $P(E,\Omega)$. 
\ed
\br
In some settings, in order to emphasize the ambient manifold $M$, we use the notation $\norm{Du}_M(\Omega)$ instead of $\norm{Du}(\Omega)$. \\
\indent For a {\Ca} set $E$ all properties are unchanged if we make alterations of any Hausdorff measure zero set. Arguing exactly as Proposition 3.1 in \cite{Giu84},  we can always choose a set $E'$ differing by a Hausdorff measure zero set with $E$ and satisfying for any $x\in \P E'$ 
\be \label{boundary}
0<vol(E'\cap B(x,\rho))<vol(B(x,\rho))
\ene 
where $\rho$ is sufficiently small. From now on,  we always assume that condition \eqref{boundary} holds for any {\Ca} set $E$. 
\er
\bd  We say a sequence of measurable functions $\{u_k\}_{k=1}^\infty$ locally converges to $u$ in $L^1(\Omega)$ if for any open set $W\subset\subset  \Omega$ we have 
\be 
\lim_{k=1}^\infty\int_W |u_k-u|dvol =0
\ene 
\ed
\br\label{rk:local:converges} If $\{u_k\}_{k=1}^\infty$ converges to $u(x)$ a.e, then this sequence locally converges to $u(x)$ in $L^1(\Omega)$. 
\er 
By the definition in \eqref{def:bv} it is easy to see that 
\bt[Lower-semicontinuity] \label{lm:bv:lower-semicontinuity}
	Let $\{u_i\}$ be a sequence of functions in $BV(\Omega)$ locally converging to a function $u$ in the $L^1(\Omega)$. Then
	\begin{equation}
		\norm{Du}(\Omega) \le \mathop {\liminf} \limits_{i\to +\infty} \norm{Du_i}(\Omega)
	\end{equation}
\et

The following $C^\infty$ approximation result is not trivial when $\Omega$ is not contained in a simply-connected domain of a Riemannian manifold. Because in this case no global symmetric mollifies exist for this domain as that in \cite[page 15]{Giu84} for Euclidean spaces. For a complete proof, we refer to \cite[section 3]{Zhou19}. 
\bt[Theorem 3.6, \cite{Zhou19}] \label{approximation_of_variation}
	Suppose $u \in BV(\Omega)$.  Then there exists a sequence of functions $\{u_k\}_{k=1}^\infty$ in $C^\infty(\Omega)$ such that
	\begin{align}
		 & \mathop{\lim} \limits_{k\to +\infty} \int_\Omega |u-u_k|d vol=0, and    \\
		 & \mathop{\lim} \limits_{k\to +\infty} \norm{Du_k}(\Omega) = \norm{Du}(\Omega)
	\end{align}
\et
As a conclusion, we can obtain the following well-known compactness result. 
\bt[Compactness theorem]\label{thm:compact:BV}
Let $\Omega \subset M$ be a bounded domain with Lipschitz boundary. Suppose there is a constant $C>0$ such that  for $\{u_k\}_{k=1}^\infty\subset BV(\Omega)$ 
\begin{equation}
\norm{u_k}_{L^1} + \norm{Du_k}(\Omega) \le C \text{ for any k}
\end{equation}
 Then there is a function $u \in BV(\Omega)$ such that there is a subsequence of $\{u_k\}_{k=1}^\infty$ converging to $u(x)$ in $L^1(\Omega)$. 
\et
We can also view a $BV_{loc}$ function as a Radon measure. 
\bt [Theorem 2.6, \cite{Zhou19}]\label{variation_measure}
	Let $\Omega \subset M$ be a bounded domain.  Suppose that $u \in BV_{loc}(\Omega)$.
	\begin{enumerate}
		\item There exists a Radon measure $ |Du|$ on $\Omega$ such that
		      \begin{equation}
			      \int_{\Omega'}f d|Du|= \sup\{\int_{\Omega'} u\div(X) dvol \mid X\in T_0\Omega' \text{ and } \la X, X\ra \le f^2\}
		      \end{equation}
		      for any open set $\Omega'\subset \subset \Omega$ and any nonnegative function $f\in L^1( |Du|)$.
		  \item There exists a vector field $\vec{\nu}$ on $\Omega$ satisfying
		      \begin{equation}
				  \int_\Omega u\div(X) dvol =-\int_\Omega \inner{X}{\vec{\nu}}d|Du|
		      \end{equation}
			  where $\la \vec{\nu},\vec{\nu}\ra=1$ $|Du|$-a.e.  for any $X \in T_0\Omega$.
	\end{enumerate}
\et
\section{The area functional}
In this section we define a new area functional $\AF(u,\Omega)$ for BV functions to generalize the area of the graph of $C^1$ functions in $M_\phi$ (see \eqref{def:area}). Suppose $N$ is a fixed Riemannian manifold with metric $\sigma$ and $\phi(x)$ is a $C^3$ positive function on $N$. Here $M_\phi$ is the set $N\PLH \R$ equipped with $\phi^2(x)(\sigma+dr^2)$. We shall establish the $C^\infty$ approximation of $\AF(u,\Omega)$ and the Miranda's observation mentioned in the introduction. \\
\indent Recall that $\Omega$ is a bounded domain in $N$, $T_0(\Omega)$ and $C_0(\Omega)$ denote the set of all smooth vector fields and smooth functions with compact supports in $\Omega$ respectively. Let $dvol$ denote the volume form of $N$. 
\bd\label{def:AF}
	Let $u$ be a measurable function in $\Omega$.  Define
	\begin{equation}
		\begin{split}
			\AF(u, \Omega)& = \sup\{\int_\Omega \phi^nh + u\div(\phi^nX)dvol \mid h \in C_0(\Omega), \\
			&X\in T_0(\Omega)\
			\text{ and } h^2+\inner{X}{X} \le 1
			\}
		\end{split}
	\end{equation}
\ed
\br
	If $u \in C^1(\Omega)$,  $\AF(u,\Omega)$ is the area of the graph of $u$ in the conformal product manifold $M_\phi$. 
	\er
\subsection{Properties of the area functional}
By the definition we have 
\bl[Lower-semicontinuity]\label{lower_semicontinuity_of_area}
Let $\{u_k\}_{k=1}^\infty$ be a sequence converging to $u(x)$ in $L^1(\Omega)$. Then
	\begin{equation}
		\AF(u,\Omega) \le \liminf_{k\to +\infty}\AF(u_k,\Omega)
	\end{equation}\el
The following lemma establishes a relation between $\AF(u,\Omega)$ and $||Du||(\Omega)$.
\bl \label{comparision_between_bv_and_area}
There is a positive constant $\mu_0:=\mu_0(\Omega, \phi)$ such that for any $u \in L^1(\Omega)$
	\begin{equation}
		\F{1}{\mu_0}\max\{\norm{Du}(\Omega), vol(\Omega)\} \le \AF(u, \Omega)\le \mu_0(\norm{Du}(\Omega)+vol(\Omega))
	\end{equation}
\el
\begin{proof}
	Since $\bar{\Omega}$ is compact, there is a positive constant $\mu_0$ such that $\F{1}{\mu_0}\le  \phi^n(x) \le \mu_0$ on $\bar{\Omega}$. By the conclusion (1) in Theorem \ref{variation_measure}  we have 
\be \label{eq:af:mid}
       \sup_{\la X,X\ra  \le 1, X\in T_0\Omega} \int_\Omega u\div(\phi^n X) dvol=\int_\Omega \phi^n(x) d|Du|
\ene  We choose $h \in C_0(\Omega)$ and $X\in T_0(\Omega)$ satisfying $h^2+\la X, X\ra \leq 1$. By the definition in \eqref{def:AF}
	\begin{equation}
		\begin{split}
			\AF(u, \Omega)
			& \le \mu_0 (vol(\Omega))+ \sup_{\norm{X}^2 \le 1} \int_\Omega u\div(\phi^n X) dvol\\
			& \le \mu_0(vol(\Omega)+||Du||(\Omega))   \text{ by \eqref{eq:af:mid} }
		\end{split}
	\end{equation}
Let $X=0$ and $h=1$, the definition \eqref{def:AF} implies that  $\AF(u, \Omega) \ge \F{1}{\mu_0}vol(\Omega)$. Now let $h=0$,  \eqref{eq:af:mid} implies that $\AF(u,\Omega)\geq \F{1}{\mu_0}||Du||(\Omega)$. The proof is complete. 
\end{proof}
 The area functional $\AF(u, \Omega)$ induces a Radon measure on $\Omega$ as follows.  
\bt \label{thm:Radon:measure} Suppose $u\in BV(\Omega)$. Then there is a Radon measure $\nu$ on $\Omega$ such that for any open set $W$ in $\Omega$, 
\be 
\nu(W)=\AF(u, W)
\ene 
\et 
\bp For any nonnegative function $f \in C_0(\Omega)$, define
\begin{equation}
	\begin{split}
		\lambda(f) & = \sup\{\int_\Omega \phi^nh + u\div(\phi^nX)dvol \mid h \in C_0(\Omega), X\in T_0(\Omega)\\
		& \text{ and } h^2+\inner{X}{X} \le f^2 \}
	\end{split}
\end{equation}
It is easy to see that $\lambda(cf)=c\lambda(f)$ for any positive constant $c>0$ and $\lambda(f+g)=\lambda(f)+\lambda(g)$.  Thus $\lambda$ is a positive linear functional on $C_0(\Omega)$.  By \cite[remark 4.3]{Sim83} (see also  \cite[Theorem 2.5]{Zhou19}), then there is a Radon measure $\nu$ on $\Omega$ such that 
\be 
\nu(W)=\sup\{\lambda(f) \mid spt(f) \subset \subset W, 0\leq f\leq 1\}
\ene 
for any open set $W\subset \Omega$. Here $spt(f)$ denotes the closure of $\{f(x)\neq 0: x\in \Omega\}$. From the definition it is clear that $\AF(u, W)=\sup\{\lambda(f) \mid spt(f)\subset \subset W, 0\leq f\leq 1\}$. The proof is complete. 
\ep 
Fix any $p\in N$. Let $exp_p$ be the exponential map at $p$. Let $inj_{\bar{\Omega}}$ be the finite number given by 
\be 
\min\{\text{ the injective radius of } x \text{ in } N: x \in \bar{\Omega}, 1\}
\ene  For any $r< inj_{\bar{\Omega}}$, let  $B_r(0)$ be the Euclidean ball centering at $0$ with radius $r$. Thus $exp_p: B_r(0)\rightarrow B_r(p)$ is a diffeomorphism. Via this exponential map we identify $B_r(p)$ with $B_r(0)$ equipped with the metric $g=g_{ij}dx^idx^j$. Then $B_r(p)$ is called \emph{a normal ball}. A vector field $X$ along $B_r(p)$ can be represented by 
\be 
    X=X^i \F{\P}{\P x_i}
\ene  
\par
Let $\vp(x)$ be a symmetric mollifier in $\R^n$, i.e. a function satisfying $\vp(x)=\vp(-x)$, $spt(\vp) \subset B_1(0)$ and $\int_{\R^n}\vp(x) dvol_{\R^n}=1$. Here $dvol_{\R^n}$ denotes $dx^1\wedge \cdots\wedge dx^n$.  If $u \in L^1(\R^n)$, let $u*\vp_\Sc$ denote the convolution of $u$ where $\vp_\Sc=\frac{1}{\Sc^n}\vp(\frac{x}{\Sc})$ namely
\begin{equation}
	u*\vp_\Sc(y)=\int_{\R^n}\vp_\Sc(x-y)u(x)dvol _{\R^n}.
\end{equation}
The convolution of a vector field $X=X^i\frac{\partial}{\partial x_i}$ and $\vp_\Sc$ on $B_r(p)$ is defined by
\begin{equation}
	X*\vp_\Sc: =X^i*\vp_\Sc\frac{\partial}{\partial x_i}
\end{equation}
Note that $dvol=\sqrt{\det(g)}dvol_{\R^n}$.  Suppose $f, h\in  C_0(B_r(p))$.   By a direct computation we have 
\begin{equation}
	\int_{B_r(p)} (\vp_\Sc*f) h dvol = \int_{B_r(p)} \frac{f}{\sqrt{\det g}} \vp_\Sc*(h\sqrt{\det g})dvol 
\end{equation}
\begin{equation}
	\int_{B_r(p)} (\vp_\Sc * h)  \div(X) dvol =\int_{B_r(p)} h \div(\frac{1}{\sqrt{\det g}} \vp_\Sc * (\sqrt{\det g}X)) dvol 
\end{equation}
To prove the $C^\infty$ approximation property of the area functional $\AF(u,\Omega)$, we need the following two techinique lemmas from \cite{Zhou19}. Both of their proofs follow from those of  \cite{Zhou19} with minor modifications. So we skip them here.
\bl[Theorem 2.12, \cite{Zhou19}] \label{uniformly_bounded}
		Let $B$ be a normal ball in $\Omega$ with a metric $g=g_{ij}dx^idx^j$.  Let $f$ be a nonnegative continuous function and $K \subset\subset  B$ be a compact subset.  Then for any $\delta >0$, there exists a $\Sc_0=\Sc_0(\Omega, K, g, f)$ such that for all $\Sc< \Sc_0$, any continuous function $h$ and any vector field $X$ satisfying
	\begin{equation}
		h^2+\la X, X\ra  \le f^2,  \text{ in } \Omega 
	\end{equation}
	we have
	\begin{equation}
		h^2+\la X', X'\ra \le (f+\delta)^2 \text{ in }  K
	\end{equation}
	where $X'$ is defined by
	\begin{equation}
		X'=\frac{1}{\sqrt{\det(g)}\phi^n(x)}\vp_\Sc*(\sqrt{\det(g)}\phi^n(x) X)
	\end{equation}
\el
\bl[Lemma 2.13, \cite{Zhou19}] \label{lm:integral}Let $B$ be a normal ball in $\Omega$ with a metric $g=g_{ij}dx^idx^j$.  Suppose $u\in BV(B)$ and $q(x)$ be a smooth function with compact support in $B$.  Then for any $\Sc>0$ there is a $\sigma_0=\sigma_0(u,g,q,K)$ such that for all $\sigma \in (0,\sigma_0)$ and any smooth vector field $X$ on $B$ satisfying $\la X, X\ra\leq 1$. 
\be 
\begin{split}
\int_B \vp_\sigma *(qu)div(\phi^n X)dvol&\leq \int_B udiv(q\phi^n Y_\sigma)dvol \\
&-\int_B \la \phi^n X, \nb q\ra dvol +\Sc
\end{split}
\ene 
where $\nb$ is the covariant derivative of $N$, $Y_\sigma=\F{1}{\sqrt{det(g)}\phi^n }\vp_\sigma*(\sqrt{det(g)}\phi^n X)$ and assume $X=0$ outside $B$. 
\el 
The $C^\infty$ approximation property of $\AF(u,\Omega)$ is stated as follows. 
\bt	\label{approximation_of_area}
	Let $\Omega$ be a bounded domain in $N$ and $u \in BV(\Omega)$.  Then there is a sequence  $\{u_k\}_{k=1}^\infty\in C^\infty(\Omega)$ such that
	\begin{align}
		 & \mathop{\lim} \limits_{k\to +\infty} \int_\Omega |u-u_k|dvol=0 \\
		 & \mathop{\lim} \limits_{k\to +\infty} \AF(u_k,\Omega)= \AF(u,\Omega)
	\end{align}
\et
\begin{proof} By Lemma \ref{comparision_between_bv_and_area} $u\in BV(\Omega)$ and that $\Omega$ is bounded imply that $\AF(u,\Omega)$ is finite. \\
	\indent By Theorem \ref{thm:Radon:measure}, there is a Radon measure $\nu$ on $\Omega$ such that $\nu(W)=\AF(u, W)$ for any open set $W\subset \Omega$. Fix any $\Sc >0$ and $r_0>0$. By  \cite[Theorem A.4]{Zhou19}, there exists a countable open cover $\{B_k\}_{k=1}^\infty$ of $\Omega$  and a positive integer $\kappa_0$ such that
	\begin{enumerate}
		\item Each $B_k$ is a normal ball such that $\nu(\partial B_k)=0$ and $\diam(B_k) \le 2 r_0$.
		\item $\{B_1,...,B_{\kappa_0}\}$ is a pairwise disjoint subcollection with
		\begin{equation}\label{eq:conclusion:one}
		\AF(u, \Omega) - \Sc \le \sum_{k=1}^{\kappa_0} \AF(u, B_k) \le \AF(u,\Omega).
		\end{equation}
		\item The subcollection $\{B_k\}_{k=\kappa_0+1}^{\infty}$ satisfies
		\begin{equation}\label{eq:conclusion:two}
		\sum_{k=\kappa_0+1}^\infty \AF (u, B_k)\le \kappa_1 \Sc
		\end{equation}
		where $\kappa_1$ is a positive integer only depending on $\Omega$ and $n$. 
		\end{enumerate}
	Take a partition of unity $\{\zeta_k\}$ subordinate to the open cover $\{B_k\}$, that is
	\begin{equation}
		\zeta_k \in C^\infty_0(B_k), 0 \le \zeta_k \le 1 \text{ and } \sum_{k=1}^{\infty}\zeta_k=1 \text{ on } \Omega.
	\end{equation}
	and for each $x\in \Omega$.  And for any $x$ in $\Omega$ there is a compact neighborhood $V$ of $x$ such that on $V$ the summation above is finite. \\
	\indent Now fix any $h\in C_0(\Omega)$ and any $X\in T_0\Omega$ satisfying $h^2+\la X,X\ra ^2\leq 1$. 	For each $k$ by Lemma \ref{uniformly_bounded} and Lemma \ref{lm:integral} we can choose $\Sc_k$ so small independent of $X$ and $h$ such that $\spt(\vp_{\Sc_k}*(u\zeta_k)) \subset B_k$,
	\begin{gather}
		\int_{B_k} |\vp_{\Sc_k} * (u \zeta_k) - u \zeta_k|d vol < \frac{\Sc}{2^k} \label{L1_convergence}  \\
(\zeta_k h)^2+\zeta_k^2\la Y_{\Sc_k}, Y_{\Sc_k}\ra \leq 1+\Sc 
\end{gather}
\be\label{eq:mid:result}
\begin{split}
\int_{B_k}\vp_{\Sc_k}*(u\zeta_k)&div(\phi^n(x)X)dvol \leq \int_{B_k} udiv(\phi^n\zeta_k Y_{\Sc_k})dvol \\
&-\int_{B_k}	u\la \phi^n X,\nb \zeta_k\ra dvol + \F{\Sc}{2^k} 
\end{split}
\ene 
Here $Y_{\Sc_k}=\frac{1}{\sqrt{\det(g)}\phi^n}\vp_{\Sc_k}*(\sqrt{\det(g)}\phi^n X)$. 
	Now we define $u_\Sc$ by
	\begin{equation}
		u_\Sc= \sum_{k}^\infty \vp_{\Sc_k} * (u\zeta_k)
	\end{equation}
	Thus $u_\Sc\in C^\infty(\Omega)$. By \eqref{eq:mid:result} we obtain 
    \begin{gather*}
		\begin{split}
			&\quad \int_\Omega \{\phi^n h + u_\Sc\div(\phi^n X)\}dvol  \\
			&=\sum_{k=1}^\infty \int_{B_k} \{\phi^n\zeta_k h+\vp_{\Sc_k} * (u\zeta_k) \div(\phi^n(x) X\} dvol \\
			&\leq \sum_{k=1}^\infty \int_{B_k}\{\phi^n\zeta_k h+udiv(\phi^n\zeta_k Y_{\Sc_k})
			-u\la \phi^n X,\nb \zeta_k\ra\} dvol+\Sc\text{ by \eqref{eq:mid:result}} \\
			&=\sum_{k=1}^\infty \int_{B_k}\{\phi^n\zeta_k h+udiv(\phi^n\zeta_k Y_{\Sc_k})\} dvol+\Sc
		\end{split}
	\end{gather*}
The last line we use the fact that $ \sum_{k=1}^\infty \int_{B_k}-u\la \phi^n X,\nb \zeta_k\ra dvol=0$.  By definition we have 
\be 
\int_{B_k}\{\phi^n\zeta_k h+udiv(\phi^n\zeta_k Y_{\Sc_k})\} dvol\leq (1+\Sc)\AF(u, B_k)
\ene 
Combining this with \eqref{eq:conclusion:one}, \eqref{eq:conclusion:two} together we obtain 
\be
\begin{split}
\int_{\Omega} \{\phi^n h + u_\Sc\div(\phi^n X)\}dvol &\leq  (1+\Sc)\sum_{k=1}^\infty \AF(u, B_k )\\
      &\leq (1+\Sc)\{\AF(u, \Omega)+\kappa \Sc \}
\end{split}
\ene 
Since $h, X$ are chosen arbitrarily satisfying $h\in C_0(\Omega), X\in T_0(\Omega)$ and $h^2+\la X, X\ra\leq 1$,  we conclude that 
\be 
\AF(u_\Sc, \Omega)\leq (1+\Sc)(\AF(u,\Omega)+\kappa \Sc)
\ene 
\indent As a result we can choose $\Sc_i\rightarrow 0$ such that $\lim_{\Sc_i\rightarrow 0}\sup\AF(u_{\Sc_i},\Omega)\leq \AF(u,\Omega)$. On the other hand by \eqref{L1_convergence} $u_{\Sc_i}$ converges to $u$ in $L^1(\Omega)$. On the other hand by Lemma \ref{lower_semicontinuity_of_area} $\lim\inf \AF(u_{\Sc_i},\Omega)\geq \AF(u,\Omega)$. Thus we obtain the conclusion. The proof is complete. 
\end{proof}
\subsection{The Miranda's observation}
In this subsection we show the Miranda's observation for $\AF(u,\Omega)$ mentioned in the introduction. \\
\indent Let $u$ be a measurable  function on $\Omega$. The set $\{(x, t)| x \in \Omega, t < u(x)\}$ is called as the \textit{subgraph} of $u$, written as $U$. For any Borel set $E$ we denote by $\lambda_E$ its characterization function. \\
\indent Let $P_\phi$ denote the perimeter of {\Ca} sets in $M_\phi$. Its connection with the area functional $\AF(u,\Omega)$ can be summarized as follows. 
\bt \label{subgraph_and_area}
Suppose $u \in BV(\Omega)$, $U$ is its subgraph and $\Omega$ is a bounded domain.  Then
\begin{equation}\label{eq:subgraph:conclusion}
P_\phi(U, \Omega\PLH\R)=\AF(u, \Omega)
\end{equation}
\et
\begin{proof} By Lemma \ref{comparision_between_bv_and_area} $\AF(u,\Omega)$ is finite.  If $u$ is a smooth function, the equality is obvious since both sides equal to the area of the graph of $u$ in $M_\phi$.  \\
	\indent Suppose $u \in BV(\Omega)$. By Theorem \ref{approximation_of_area} there is  a sequence of smooth functions $\{u_k\}_{k=1}^\infty$ such that	$\{u_k\}$ converges to $u$ \text{ in } $L^1(\Omega)$ and $\AF(u_k, \Omega)$ converges to $ \AF(u, \Omega)$ as $k\rightarrow \infty$.  Let $U_k$ be the subgraph of $u_k$.  Note that $\{\lambda_{U_k}\}$ converges  to $\lambda_{U}$ in a.e.  By Lemma \ref{lm:bv:lower-semicontinuity} we have
	\be \label{eq:sdt}
	\begin{split}
P_\phi(U,\Omega\PLH\R) & \le \lim_{k\rightarrow \infty}\inf P_\phi (U_k,\Omega \PLH \R)\\
	&=\lim_{k\rightarrow \infty}\inf\AF(u_k, \Omega) = \AF(u, \Omega)
	\end{split}	
\ene 
\indent To prove the converse first suppose that $|u| \le \mu$ for some positive constant $\mu>1$.  Observe that we can add a constant to $u$ without changing $P_\phi(U,\Omega \PLH \R)$ and $\AF(u, \Omega)$. Thus without loss of generality we assume that $u \ge 1$.  Let $X \in T_0\Omega$ and $h(x)\in C_0(\Omega)$ satisfying $\la X, X\ra + h^2(x) \le 1$.  Let $\eta(t)$ be a smooth nonnegative function with compact support in $(0, \mu +1)$ such that $\eta \equiv 1$ in $[1, \mu ]$ and $\eta \le 1$ on $(0,\mu+1)$.  Define
	\begin{equation}
	X'=\frac{1}{\phi(x)}\eta(r)(X+h(x)\partial_r)
	\end{equation}
	Then $X' \in T_0(\Omega \PLH \R)$ and $\la X', X'\ra \le 1$ in $M_\phi$.  Observe that for each $x\in \Omega$, we have
	\be
	\int_\R \lambda_U(x,r) \eta dr=u(x)+ C ,\quad  \int_\R \lambda_U(x,r) \eta' dr=1
\ene 
	where $C$ is a constant depending only on $\eta$.  We denote the product manifold $N\PLH\R$ with metric $\sigma+dr^2$ by $M$. Let $div_{M_\phi}$ and $div_M$ be the divergence on $M_\phi$ and $M$ respectively.  Note that $div_{M_\phi}(X')=\frac{1}{\phi^{n+1}}div_M{(\phi^{n+1}X')}$ (For example see \cite[(2.11)]{Zhou19}.  Then by definition of perimeter, we have 
\be\notag
	\begin{split}
	P_\phi(U, \Omega \PLH \R) &\ge \int_{\Omega \PLH \R}\lambda_U \div_{M_\phi}(X')dvol_{M_\phi}\\
	&=\int_\Omega \{ \div( \phi^n X) \int_\R\lambda_U \eta dr + h\phi^n (\int_\R \lambda_U \eta' dr )\}dvol  \\
	&=\int_\Omega \{u \div(\phi^n X) + h \phi^n \}dvol 
	\end{split}
	\ene 
	Taking supremum over both sides for all $h\in C_0(\Omega), X\in T_0\Omega$ with $h^2+\la X, X\ra\leq 1$ yields that $P_\phi(U, \Omega\PLH\R)\ge \AF(u, \Omega)$. With \eqref{eq:sdt} this yields \eqref{eq:subgraph:conclusion} in the case of $|u|\leq \mu$. \\
	\indent  As for the general case we consider an approximation procedure based on the finite case. Let $u_T$ be the truncation of $u$ by $T$, i.e. $u_T=\max\{\min\{u, T\}, -T\}$.  Let $U_T$ be the subgraph of $u_T$. Thus $\lambda_{U_T}$ converges to $\lambda_{U}$ a.e. in $\Omega\PLH \R$. By \eqref{eq:sdt} and the lower semicontinuity in Lemma \ref{lower_semicontinuity_of_area}, $P_\phi(U, \Omega\PLH \R)$ is finite. Note that 
	\be
	\begin{split} \lim_{T\rightarrow +\infty}&P_\phi (U,  \Omega\PLH (-T, T))=P_\phi (U, \Omega\PLH \R)\\
	 &P_\phi (U_T, \Omega\PLH(-T, T))=P_\phi (U, \Omega\PLH(-T,T))
	 \end{split}
	\ene  
	Combining this with Lemma \ref{lower_semicontinuity_of_area} we obtain 
	\be\label{eq:mist:uv}
	\begin{split}
\lim_{T\rightarrow +\infty }P_\phi(U_T, \Omega\PLH \R )=P_\phi(U, \Omega\PLH(-T,T))
\end{split}
\ene
As a result we conclude 
	\begin{equation}
	\begin{split}
	P_\phi(U,\Omega\PLH\R)&=\lim_{T\rightarrow +\infty}P_\phi(U_T, \Omega\PLH\R)\\
	&=\lim_{T\rightarrow +\infty} \AF(u_T, \Omega) \ge \AF(u, \Omega)
	\end{split}
	\end{equation}
	Combining this with \eqref{eq:sdt} yields the conclusion. The proof is complete. 
\end{proof}
As an application we obtain a decreasing property of certain {\Ca} sets in $M_\phi$. 
\bt  \label{thm:decrease_perimeter}
	Let $E \subset\Omega\PLH \R$ be a Caccioppoli set in $M_\phi$ with the following assumption: for almost every $x\in \Omega$, there exists a $T_x>0$ such that $\lambda_E(x,t)=0$ for all $t>T_x$ and $\lambda_E(x,t)=1$ for all $t<-T_x$. Then the function
	\begin{equation}\label{def:w}
		w(x)=\lim_{k\to +\infty} (\int^k_{-k} \lambda_E(x, t)dt -k)
	\end{equation}
is well-defined and 
	\begin{equation}\label {eq:uset:rt}		
	\AF(w, \Omega) \le P_\phi(E, \Omega\PLH \R)
	\end{equation}
\et
\br\label{Rm:decreasing:perimeter} This result should be firstly observed by Miranda \cite{Mir64} in the case of $\phi(x)\equiv 1$ and $N=\R^{n}$ (see also \cite{Alm87}). It is a generalization  of   \cite[Theorem 14.8]{Giu84} and the symmetrization of hyperbolic spaces in  \cite[remark 2.3]{Lin85} by Lin. It is also similar to the arrangement of constant mean curvature functional in \cite[Theorem 3.1]{DS09} by De Silva-Spruck, the decreasing perimeter property of singular area functional  in \cite[Lemma 9]{BD87} by Bemelmans-Dierkes, \cite[Lemma 3.3]{Ger98} by Gerhardt and  \cite[Lemma 5.8]{Zhou19}. 
\er 
\begin{proof}  By the assumption on $E$ it is not hard to see that $\omega(x)$ is well defined. Moreover we assume that $E$ has finite perimeter in $M_\phi$. Otherwise nothing needs to prove. \\
	\indent First we assume that there is a $T>0$ such that $\P E\subset \Omega\PLH(-T, T)$.  Let $\eta(r)$ be a compactly supported smooth function on $\R$ such that $\eta(t)\equiv 1$ on $(-1,1)$ and $\eta(r) \le 1$.  For each $k \in N^+$, define
	\begin{equation}
		\eta_k(r)= \left\{
		\begin{aligned}
			 & \eta(r+k) \quad r\in (-\infty, -k) \\
			 & 1 \quad r\in [-k, k]               \\
			 & \eta(r-k) \quad r\in (k, +\infty)
		\end{aligned}
		\right.
	\end{equation}
	Choose $h \in C_0(\Omega)$ and $X \in T_0(\Omega)$ satisfying $h^2+\la X, X\ra \le 1$. Define 
	\begin{equation}
		X'=\frac{1}{\phi(x)}\eta_k(r)(X+h(x)\partial_r)
	\end{equation}
Thus $\la X', X'\ra \le 1$.  Then we have
	\begin{equation}
		\begin{split}
			&P_\phi(E, \Omega\PLH \R) \ge \int_{\Omega \PLH \R} \lambda_E \div_{M_\phi}(X')dvol_{M_\phi}  \\
			&=\int_{\Omega \PLH \R} \lambda_E \frac{1}{\phi^{n+1}}div_M(\phi^{n}\eta_k(r)(X+h\partial_r))\phi^{n+1}dr dvol \\
			&=\int_\Omega \{\div(\phi^n X) \int_{\R}\lambda_E \eta_k(r)dr+ \phi^n h \int_\R \lambda_E \eta_k'(r)dr\} d vol 
		\end{split}
	\end{equation}
	where $dvol_{M_\phi}$ and $div_{M_\phi}$ denote the volume form and the divergence of $M_\phi$ respectively, $div_M$ is the divergence of the product manifold $M:=(N\PLH\R, \sigma+dr^2)$.  By the assumption on $E$ and $\P E\subset \Omega \PLH [-T, T]$, for almost every $x \in \Omega$ we have
	\begin{align}
		&\lim_{k\to +\infty}\int_{\R}\lambda_E(x,r)\eta'_k(r)dr = 1 \\
		&\lim_{k\to +\infty}\int_{\R} \eta_k(r) \lambda_E(x, r)dr= w(x)+ C 
	\end{align}
	where $C$ is a constant depending only on $\eta$ and $T$. \par
	Let $k \to \infty$, we get
	\begin{equation}
		P_\phi(E, \Omega \PLH \R) \ge \int_\Omega \div(\phi^n X)w+\phi^n h d vol 
	\end{equation} 
	Because $h(x)\in C_0(\Omega)$ and $X\in T_0\Omega$ are arbitrarily choosen satisfying $h^2+\la X, X\ra \leq 1$, we obtain \eqref{eq:uset:rt} in the case of $\P E\subset \Omega\PLH(-T, T)$. \\
	\indent As for a general {\Ca} set $E$ with finite perimeter, we set $E_T=E\cup \Omega\PLH(-\infty,-T)\backslash\Omega \PLH [T, +\infty)$. Note that $\lim_{T\rightarrow +\infty}P_\phi (E, \Omega\PLH (-T, T))=P_\phi (E, \Omega\PLH \R)$ and $P_\phi (E_T, \Omega\PLH(-T, T))=P_\phi (E, \Omega\PLH(-T,T))$. By the lower semicontinuity in Lemma \ref{lower_semicontinuity_of_area}, we obtain
	\be 
	\lim_{T\rightarrow +\infty} P_\phi(E_T, \Omega\PLH \R)=P_\phi(E, \Omega\PLH \R)
	\ene 
	 Define $$w_T(x)=\lim_{k\to +\infty} (\int^k_{-k} \lambda_{E_T}(x, t)dt -k)$$
	 By the assumption of $E$, $w_T(x)$ converges to $w(x)$ a.e. $x\in\Omega$ as $T\rightarrow +\infty $. Again by Lemma \ref{lower_semicontinuity_of_area} 
	 \be 
	 \begin{split}
	   \AF(w(x), \Omega) &\leq\lim\inf\AF(w_T(x), \Omega)\\
	                    &\leq \lim\inf P_\phi(E_T,\Omega\PLH \R)\leq P_\phi ( E,\Omega\PLH\R)
	   \end{split}
	 \ene 
	 The proof is complete. 
\end{proof}
Theorem \ref{thm:decrease_perimeter} and the following result together are referred as the Miranda's observation.
\bt Let $u \in BV(\Omega)$ with $\AF(u,\Omega)<\infty$. Suppose $u(x)$ locally minimizes $\AF(.,\Omega)$, i.e. $\AF(u, \Omega)\leq \AF(v,\Omega)$ if  $U\Delta V\subset \subset\Omega \PLH\R$ where $U,V$ are the subgraphs of $u(x), v(x)$ respectively. Then $U$ locally minimizes the perimeter in  $M_\phi$, i.e.
          $$P_\phi(U, \Omega \PLH\R)\leq P_\phi(F, \Omega \PLH\R)$$ 
        for any {\Ca} set $F$ satisfying $F\Delta U \subset\subset \Omega \PLH\R$. 
\et 
\bp
	Let $F \subset \Omega \PLH \R$ be a Caccioppoli set satisfying $F\Delta U\subset \subset \Omega\PLH \R$. Since $U$ is a subgraph of $u(x)$, it is easy to see $F$ satisfying the condition in Theorem \ref{thm:decrease_perimeter}. Let $w(x)$ be defined as in Theorem \ref{thm:decrease_perimeter}.  Then
	\begin{equation}
		P_\phi(U, \Omega \PLH \R) = \AF(u, \Omega) \le \AF(w(x), \Omega) \le P_\phi(F, \Omega \PLH \R)
	\end{equation}
	The poof is complete. 
\ep
\section{The area minimizing problem}
In this section we consider the area minimizing problem in the conformal product manifold $M_\phi$ (see Definition \ref{eq:Model}) similar to that in \cite{GZ20}. \\
\indent Throughout this section suppose $\Omega$ and $\Omega'$ are two bounded $C^2$ domains in $N$ satisfying $\Omega\subset \subset \Omega'$. Let $Q_\phi$ be the set $\Omega\PLH\R$ in  $M_\phi$ and $\bar{Q}_\phi$ be its closure $\bar{\Omega}\PLH\R$.
 Suppose $\psi(x)\in C^1(\Omega'\backslash\Omega)$. We denote by $n$ the dimension of $N$. Let $\mG$ denote the set of $n$-{\IMC}s with compact support in $\bar{Q}_\phi$, i.e. for any $T\in \mG$, its support $spt(T)$ is contained in $\bar{\Omega}\PLH [a,b]$ for some finite numbers $a<b$.  Let $\Gamma$ be the graph of $\psi(x)$ on $\P\Omega\PLH\R$.   The area minimizing problem in this section is to find a solution to attain the value 
\be \label{mass:minimizing}
\min\{\bM(T): T\in  \mG, \P T=\Gamma \}
\ene
where $\bM$ is the mass of $T$ in $M_\phi$ (see subsection 4.1). \\
\indent The main result of this section is stated as follows. 
\bt\label{thm:area:minimizing} Suppose $\phi(x)$ is a $C^3$ positive function on $N$. Then there is a $u(x)\in BV(\Omega')$ satisfying $u(x)=\psi(x)$ outside $\Omega$ with the following three properties:
\begin{enumerate}
	\item $u(x)$ realizes  the minimum 
\be\label{functional:area:minimizing}
\min \{\AF(v(x),\Omega'): v(x)\in  BV(\Omega'), v(x)=\psi(x) \text{  outside $\Omega$}\}
\ene
\item  $T=\P [[U]]|_{\bar{Q}_\phi}$ solves the area minimizing problem \eqref{mass:minimizing} where $U$ is the subgraph of $u(x)$ in $\Omega'\PLH\R$,
\item  $u(x)\in C^{2,\beta}(\Omega)$ satisfies \eqref{eq:minimal:surface:equation:phi} and $\min_{\P \Omega}\psi(x)\leq u(x)\leq  \max_{\P \Omega}\psi(x)$ for any $x\in \Omega$ and any $\beta\in (0,1)$. 
\end{enumerate}
  \et  
  \br In the case of $\phi(x)\equiv 1$, the above theorem is obtained by Lin in \cite[section 4.1]{Lin85}. But we do not understand the proof of its uniqueness just assuming $\Omega$ is bounded Lipschitz. 
  \er 
  \br \label{Rm:boundary} We show that for $T$ in the property (2) $\P T=\Gamma$. For $r>0$ we denote the set $\{x\in \Omega': dist(x,\Omega)<r\}$ by $\Omega_r$. When $r$ is sufficiently small, $\Omega_r$ is a $C^2$ domain.  We define $T_r:=\P [[U]]|_{\bar{\Omega}_r\PLH\R}$. It is obvious that $T_r$ is a $C^2$ graph near $\P \Omega_r$. Thus $\P  T_r = graph(\psi(x)|_{\P\Omega_r})$. And $T_r$ converges to $T$ and $\P  T_r$ converges to $\P T$ in the sense of current as $r\rightarrow 0$. Moreover $graph(\psi(x)|_{\P\Omega_r})$ converges to $\Gamma$ in the $C^1$ sense as $r\rightarrow 0$. By the definition of the current $\P T=\Gamma$. This is
  the reason that we assume $\psi(x)\in C^1(\Omega'\backslash \Omega)$. 
  \er 
\subsection{Integer Multiplicity Currents}\label{sec:IMC} Here we collect some necessary facts on {\IMC}s. Our main references are \cite[section 27]{Sim83} by Simon and  \cite{LY02} by Lin-Yang. \\
\indent  Let $U$ be an open domain in a Riemannian manifold $M$ with dimension $m$ and $H^j$ denotes the $j$-dimensional Hausdorff measure in $M$ for any $j>0$. Suppose $k$ is an integer in $[0, m]$. Let $D^k(U)$ be  the set of all $k$-smooth forms with compact support in $U$. A $k$-current in $M$ is a linear functional in $D^{k}(M)$. 
\begin{Def} A set $E\subset M$ is said to be countable $k$-rectifiable if 
	$$ E\subset E_0\cup_{j=1}^\infty F_j(E_j)$$
	where $H^k(E_0)=0$ and $F_j:E_j\subset \R^k\rightarrow M$ is a Lipschitz map for each $j$. 
\end{Def}
Now we can define a $k$-integer multiplicity rectifiable current. 
\begin{Def} Let $T$ be a $k$-current in $M$. We say that $T$ is a $k$-integer multiplicity current if 
	\be 
	T(\omega)=\int_S\la \omega,\eta\ra \theta(x)d H^k(x) 	\ene 
	where $S$ is a countable $k$-rectifiable subset of $M$, $\theta$ is a positive locally $H^k$-integrable function which is integer-valued, and $\eta$ is a $k$-form $\tau_1\wedge\cdots\wedge \tau_k$ oriented the tangent space of $S$ a.e. $H^k$. $T$ is also written as
	$\tau(S,\theta, \eta)$. \\
	\indent The mass of $T$ in $U$ is 
	\be
	\bM_U(T):=sup\{T(\omega): \la \omega, \omega\ra \leq 1, \omega \in D^k(U) \}
	\ene 
	where $<,>$ denotes the usual pairing of $k$-form. The boundary of $T$ is defined by $\P T(\omega)=T(d\omega)$ for any $\omega \in D^{k-1}(U)$. 
\end{Def}
\br \label{rk:mark:notation} For any $k$-submanifold $M'$ $[[M']]$ is a k-{\IMC} just choosing $\eta$ as its orientation which is equal to $\tau(M', 1,\eta)$\\
\indent If the dimension of an {\IMC} $T$, written as $\tau(S, \theta,\eta)$, is equal to the dimension of $M$ we always choose $\eta$ as the volume form of $M$. In this case $T$ is written as $\tau(S,\theta)$. 
\er 
A good property of {\IMC}s is their compactness obtained by Federer and Fleming in \cite{FF60} (see also \cite[Theorem 27.3]{Sim83}).
\bt\label{compact:thm} Suppose $\{T_j\}_{j=1}^\infty$ is a sequence of $k$-{\IMC}s with $$\sup\{\bM_W(T_j)+\bM_W(\P T_j)\}<\infty$$ for any open set $W\subset\subset M$. Then there is a $k$-{\IMC} $T$ such that $T_j$ converges weakly to $T$ and 
$$
\bM_W(T)\leq \lim_{j\rightarrow +\infty} \sup_{i\geq j} \bM_W(T_i)
$$
 \et 
A useful way to construct {\IMC}s is the pushforward of local Lipschitz maps. 
\begin{Def} Let $U, V$ be two open sets in (different) Riemannian manifolds. Suppose $f:U\rightarrow V$ is local Lipschitz,  $T=\tau(S, \eta, \theta)$ is a $k$-{\IMC} and $f|spt T$ is proper.  We define $f_\#T$ by
	\be 
	f_{\#}T(\omega)=\int_S\la \omega|_{f(x)}, df_{\#}\eta\ra \theta(x) dH^k(x)
	\ene 
\end{Def}
    \bt \label{thm:closure} Let $Q_\phi$ be the conformal cone in $M_\phi$ (see Definition \ref{def:Model}) and $n=dim N$. Let $k<n+1$ be a positive integer. Let $T$ be a $k$-{\IMC} in $\mathcal{G}$ satisfying $\P T =0$. Then there is a $(k+1)$-{\IMC} $R$ in $\bar{Q}_\phi$ such that $\P R =T$. Here $spt(R)$ may be noncompact in $\bar{Q}_\phi$. 
\et 
\br
A similar proof also appears in \cite[section 26.26]{Sim83}.   
\er
\bp The proof is exactly the same as that of  \cite[Theorem 3.9]{GZ20}. Note that even with the metric $\phi^2(x)(\sigma+dr^2)$ the map $h:M_\phi \PLH \R \rightarrow M_\phi$ as $h((x,r),t)=(x,r+t)$ is still proper and local Lipschitz. Thus $h_{\#}([[(-\infty,0)]]\PLH T)$ is well-defined. Moreover 
\begin{align*}
\P h_{\#}([[(-\infty, 0)]]\PLH T)&=h_{\#}(\P ([[(-\infty, 0)]]\PLH T))\\
&= h_{\#}(\{0\}\PLH T)-h_{\#}((-\infty,0)\PLH \P T)\\
&=T-h_{\#}((-\infty,0 )\PLH \P T)
\end{align*}
The conclusion follows from that $\P T=0$. 
\ep 
 Fix two numbers $a<b$ such that $\Gamma \subset \P\Omega\PLH (a,b)$. Now consider an auxiliary problem of \eqref{mass:minimizing} as follows. 
\be \label{au:mass:minimizing}
A_{a,b}:=\min\{\bM(T): T\in  \mG, \P T=\Gamma, spt(T)\subset \bar{\Omega}\PLH [a,b] \}
\ene
By Theorem \ref{compact:thm} there is a n-{\IMC} $T_{a,b}$ contained in  $\bar{\Omega}\PLH[a,b]$ with $\P T_{a,b}=\Gamma$ satisfying $\bM(T_{a,b})=A_{a,b}$. We denote $\Omega'$ by an open domain in $N$ satisfying $\Omega\subset\subset \Omega'$. Similar to  \cite[Lemma 4.6]{GZ20} we have 
\bl\label{lm:au:mass:minimizing} We can choose $T_{a,b}$ such that there is a {\Ca} set $F$ in $\Omega'\PLH \R $ such that $T_{a,b}=\P [[F]]|_{\bar{Q}_\phi}$  and $F$ is the subgraph of $\psi(x)$ outside $\bar{Q}_\phi$ .   \el 
\br\label{Rm:key:fact}  Because  $T_{a,b}=\P [[F]]|_{\bar{Q}_\phi}$ is contained in $\bar{\Omega}\PLH [a,b]$ one sees that for each $x\in \Omega$ it holds that $\lambda_{F}(x, t)=1$ for $t<a$ and $\lambda_{F}(x, t)=0$ for $t>b$. If not the case we can replace $F\cap \Omega\PLH\R$, the complement of $F$, with $F^c\cap \Omega\PLH \R$ without any change of the perimeter. 
\er 
\bp  Since $\psi(x)\in C^1(\Omega'\backslash\Omega)$ we can extend $\psi(x)$ as a $C^1$ function on $\Omega'$ satisfying $a<\psi(x)<b$ on $\Omega'$ and its subgraph $E$ is a {\Ca}  set in $\Omega'\PLH\R$. Define $S:=\P [[E]]|_{\bar{\Omega}\PLH[a,b]}$ . Then by remark \ref{Rm:boundary} $\P S=\Gamma$ .\\
\indent From Theorem \ref{thm:closure} there is a $(n+1)$-{\IMC} $R$ in $\Omega'\PLH (-\infty,b]$ such that 
$T_{a,b}-S=\P R$. Then we have 
\be\label{eq:reason:first}
T_{a,b}=\P [[E]]|_{\bar{Q}_\phi}+\P R=\P [[E]]|_{\bar{\Omega}\PLH[a,b]}+\P R|_{\bar{\Omega}\PLH[a,b]}
\ene 
Observe that $[[E]]+R$ can be represented as $\tau(\Omega'\PLH (-\infty,b), \theta)$ where $\theta$ is some integer value measurable function on $N\PLH \R$. Since $spt(R)\subset \bar{\Omega}\PLH (-\infty, b]$, $\theta =\lambda_E$ outside $\bar{Q}_\phi$ i.e.
\be \label{evalue:theta}
\theta =1 \text{ or }  0 \quad \text{ outside $\bar{Q}_\phi$}
\ene Now define a function $\theta_1=\theta$ if $\theta\neq 1$ and $\theta_1=0$ if $\theta=1$. Let $\theta_0$ denote the function of $\theta-\theta_1$. \\
\indent Set $F:=\{p\in \Omega'\PLH(-\infty,b]: \theta(p)=1\}$. It is not hard to see that $\theta_0=\lambda_F$ and 
\be\label{det:F}
[[F]]=\tau(\Omega'\PLH(-\infty, b), \theta_0)
\ene  As a result one has the following decomposition 
\be \label{eq:reason:two}
[[E]]+R=[[F]]+G
\ene 
where $G$ is the {\IMC} $\tau(\Omega'\PLH (-\infty, b),\theta_1)$. By \eqref{evalue:theta} and the definition of $\theta_1$, $\theta_1=0$ outside  $\bar{Q}_{\phi}$. As a result $spt( G)\subset \bar{Q}_\phi$. \\
\indent Now define  $U_j=\{ p\in \Omega'\PLH (-\infty,b): \theta(p)\geq j\}$ for any integer $j$. By the definition of $E$, we have $spt(\P[[U_j]])\subset \bar{Q}_\phi$ for any $j\neq 1$. Note that $spt(T_{a,b})\subset\bar{\Omega}\PLH [a,b]\subset  \bar{Q}_\phi $.  Applying the decomposition theorem of {\IMC}s \cite[corollary 27.8]{Sim83} we obtain 
\begin{align*}
\mu_{ T_{a,b}}&=\sum_{j=-\infty,j\neq 1}^\infty\mu_{\P [[U_j]]}+\mu_{\P  [[U_1]]}|_{\bar{Q}_\phi}\\
&=\sum_{j=-\infty,j\neq 1}^\infty\mu_{\P[[ U_j]]}|_{\bar{\Omega}\PLH[a,b]}+\mu_{\P  [[U_1]]}|_{\bar{\Omega}\PLH[a,b]}
\end{align*}
This implies that 
\be\label{eq:property:G}
spt(\mu_{\P [[U_j]]})\subset \bar{\Omega}\PLH [a,b], j\neq 1, spt(\mu_{\P[[ U_1]]})|_{\bar{Q}_\phi}\subset  \bar{\Omega}\PLH [a,b]
\ene 
For $G=\tau(M,\theta_1)$ applying the decomposition theorem \cite[corollary 27.8]{Sim83} again gives that 
\be 
\mu_{\P G}=\sum_{j=3}^\infty \mu_{\P[[ U_j]]}+2\mu_{\P [[U_2]]}+\sum_{j=-\infty}^0 \mu_{\P [[U_j]]},  
\ene 
\indent 
By \eqref{eq:property:G} $spt (\P G)\subset \bar{\Omega}\PLH[a,b]$. As for $F$ the decomposition theorem gives that $\P [[F]]=\P [[U_1\backslash U_2]]$. Since $U_2\subset U_1$, with \eqref{eq:reason:first} and \eqref{eq:reason:two},  \cite[Lemma 3.17]{GZ20} implies that 
\be\label{eq:fact:F}
\mu_{\P [[F]]}\leq \mu_{\P [[U_1]]}+ \mu_{\P [[U_2]]}   \ene 
and $spt(\P [[F]])|_{\bar{Q}_\phi}\subset \bar{\Omega}\PLH [a,b]$ from \eqref{eq:property:G}.  Thus $\bM(\P [[F]]|_{\bar{Q}_\phi})\leq  \bM(T_{a,b})<\infty$. From \eqref{eq:reason:two}, 
\be 
T_{a,b}= \P [[F]]|_{\bar{Q}_\phi}+\P G=\P [[F]]|_{\bar{\Omega}\PLH[a,b]}+\P G
\ene 
and $\P (\P [[F]]|_{\bar{Q}_\phi})=\P T_{a,b}=\Gamma$.  Since $T_{a,b}$ solves \eqref{mass:minimizing}, \eqref{eq:fact:F} implies that we can choose $T_{a,b}$ as $\P [[F]]|_{\bar{Q}_\phi}$. By \eqref{eq:reason:two} and \eqref{det:F}  $F$ coincides with the subgraph of $\psi(x)$ outside $\bar{Q}_\phi$. 
The proof is complete. 
\ep 
\subsection{Regularity of almost minimal boundary}
In this subsection we recall some results on the regularity of the almost minimal boundary from \cite[section 3.2]{GZ20}.  All of their proofs are skipped here. We refer the reader to  \cite{Giu84}, \cite{Tam82} and \cite{DS93} for more details. Throughout this subsection let $M$ be a $n+1$-dimensional Riemannian manifold. \\
\indent Suppose $T$ is a $(n+1)$-dimensional  {\IMC}  in $M$, represented as  $\tau(V,\theta)$ where $V$ is a $L^{n+1}$ measurable set of $M$. The relationship between the mass and BV functions can be summarized as follows:
\be\label{Mass_and_min}
\bM_W(\P T)=||D\theta||_{M}(W);  
\ene 
for any open set $W\subset\subset M$.  In particular if $E$ is a {\Ca} set in $W$ and $T=[[E]]$  we have 
\be \label{eq:relation:st}
\bM_W(\P  [[E]])=||D\lambda_E||_M(W)=P(E, W)
\ene 
where $P$ is the perimeter in $M$. For a derivation see \cite[subsection 27.7]{Sim83}. \\
\indent For a point $p\in M$ and $r>0$ we denote by $B_r(p)$ the open ball centered at $p$ with radius $r$ . Now we define an almost minimal set in an open set and a closed set respectively. 
\begin{Def}\label{key:def}
Let $\Omega$ be a domain. Fix $\mu_0>0$ and $\alpha\in (0,\F{1}{2}]$.	Suppose $E\subset \Omega$ is a {\Ca} set in $M$. 
	\begin{enumerate} 
		\item We say that $E$ is an $(\mu_0,\alpha)$-almost minimal set in $\Omega$ if there is an $r_0>0$ and a constant $\mu_0$ with the property that for any $r< r_0$, $x \in \Omega$ with $B_r(p)\subset\subset \Omega$, 
		\be \label{DEF:a}
		P(E, B_r(x))\leq P(F,B_r(x))+\mu_0r^{n+2\alpha}
		\ene 
		where $F$ is any {\Ca} set satisfying $E\Delta F\subset B_r(p)$. In particular if $C=0$, we say $F$ is a minimal set in $\Omega$. 
		\item 	We say that $E$ is an $(\mu_0,\alpha)$-almost minimal set in $\bar{\Omega}$ if there is an $r_0>0$ and a constant $\mu_0$ with the property that for any $r< r_0$, $x \in \bar{\Omega}$ 
		\be 
		P(E, B_r(x))\leq ||D\lambda_{F\cap\Omega}||_M(B_r(x))+\mu_0r^{n+2\alpha}
		\ene 
		where $F$ is any {\Ca} set satisfying $E\Delta F\subset B_r(p)$. In particular if $C=0$, we say $F$ is a minimal set in $\Omega$. 
		\item  The regular set of $\P  E$ is the set $\{p\in \P E: \P E$ is a $C^{1,\alpha}$ graph in a ball containing $p \}$.  The singular set of $\P  E$ is the complement of the regular set in $\P E$. 
	\end{enumerate}
\end{Def}
\br  \label{rk:ctwo}By \cite[Lemma 7.6]{Zhou19} all $C^2$ bounded domains are almost minimal sets in an open neighborhood of their boundaries. 
\er 
A good property of almost minimal sets in a domain is their boundary regularity. 
\bt[Theorem 1 in \cite{Tam82}, Theorem 5.6 in \cite{DS93}] \label{thm:regularity:minimal:set}Suppose a {\Ca} set $E$ is a $(\mu_0,\alpha)$-almost minimal set in a domain $\Omega$. Let $S$ be the singular set of $\P  E$ in $\Omega$. Then 
\begin{enumerate}
	\item  if $m\leq 7$, $S=\emptyset$;
	\item if $m=7$, $S$ consists of isolated points;
	\item if $m>7$, $H^{t}(S)=0$ for any $t>m-7$.  Here $H$ denotes the Hausdorff measure in $M$. 
\end{enumerate}
where $m=dim M$. 
\label{regularity:accounting:thm}
\et 
\br Note that in general the boundary of almost minimal sets in closed sets does not have such good regularity. 
\er 
The following result is very important. 
\bt [Theorem 3.19 in \cite{GZ20}] \label{regularity:key:thm}  Let $\Omega_1$ and $\Omega_2$ be two $C^2$ domains in $M$. Define $\Omega'=\Omega_1\cap \Omega_2$.  Fix a point $p$ in $\P\Omega_1\cap \P \Omega_2$.  Suppose $E\subset \Omega'$ is an $(\mu_0,\alpha_0)$- almost minimal set in $\bar{\Omega}'$ (the closure of $\Omega'$) and $\P E$ passes through $p$, then $\P E$ is a $C^{1,\alpha'}$ graph in an open ball containing $p$ for some $\alpha'\in (0, 1)$. 
\et 
\br  In general the boundary of $\Omega'$ is not $C^2$.
\er 
 
    \subsection{The proof of Theorem \ref{thm:area:minimizing}}  
      \bp  Set $\alpha_0:= \min_{\P \Omega}\psi(x)$ and $\alpha_1: =\max_{\P\Omega}\psi(x)$.  Fix $n>0$ such that $-n<\alpha_0\leq \alpha_1<n$. 
      By Lemma \ref{lm:au:mass:minimizing} there is a {\Ca} set in $\Omega'\PLH \R$ such that $T_{n}=\P [[F_n]]|_{\bar{\Omega}\PLH [-n,n]}$ realizes the minimum in \eqref{au:mass:minimizing} and $F_n$ is the subgraph of $\psi(x)$ outside $\bar{Q}_\phi$. \\
      \indent We claim that $T_{n}$ has to be contained in the closed set $\bar{\Omega} \PLH [\alpha_0, \alpha_1]$. If not the case, assume 
      \be 
      r_0:=\max\{r\in \bar{\Omega}: p=(x,r)\in T_{n}=\P [[F_n]]|_{\bar{Q}_\phi}\}\in (\alpha_1, n]
      \ene 
      Let $p_0$ be the point on $T_n$ achieving this maximum. Since $r_0>\alpha_1$, there is an embedded ball $B_r(p_0)$ such that $B_r(p_0)\cap \Gamma=\emptyset$.  \\
      \indent The first case is that $p_0=(x_0, r_0)$ for some $x_0\in \P\Omega$. Note that $F_n$ is an minimal set in $\bar{\Omega}\PLH \R\cap B_{r}(p_0)$ By Theorem \ref{regularity:key:thm} $\P F_n$ is still a $C^{1,\alpha'}$ graph near $p_0$ for some $\alpha'\in (0,1)$. Moreover $\P F_n$ is tangent to $\P\Omega\PLH \R$ and $\Omega \PLH \{r_0\}$ at $p_0$. But $\P \Omega \PLH \R$ and $\Omega\PLH \{r_0\}$ is orthogonally transverse at $p_0$. This gives a contradiction and the first case is impossible.  \\
      \indent The second case is $p_0=(x_0, r_0)$ for some $x_0\in \Omega$. By remark \ref{Rm:key:fact} $\P [[F_n]]|_{\bar{Q}_\phi}\subset \bar{\Omega}\PLH [-n,n]$. Thus $F_n$ is a minimal set in $\bar{\Omega}\PLH [-n,n]$ away from $\Gamma$. Note that  $\Omega\PLH[-n,n]$ be the intersection of two $C^2$ domain.  Again by Theorem \ref{regularity:key:thm} $\P F_n$ is a $C^{1,\alpha'}$ graph in a neighborhood of $p_0$ contained in $B_r(p_0)$.  But we should observe that $\Omega'\PLH \{r_0\}$ is also minimal in $M_\phi$ by Theorem \ref{thm:minimal:equation}. With respect to the upward normal vector, $H_{\P F_n}\leq 0$ near $p_0$ in the Lipschitz sense. By the maximum principle in  \cite[Theorem A.1]{GZ20} $\P F_n$ coincides with $\Omega' \PLH \{r_0\}$ near $p_0$. From the connectedness of $\Omega\PLH\{r_0\}$ we obtain that $\Omega'\PLH \{r_0\}\subset \P F_n$. This is impossible since $\P F_n$ is contained in $\bar{\Omega}\PLH[-n,n]$. Thus the second case is also impossible. This means $r_0\leq \alpha_1$. Arguing a similar derivation we will obtain that 
      \be
      \min\{r\in \bar{\Omega}: p=(x,r)\in T_{n}=\P [[F]]|_{\bar{Q}_\phi}\}\geq \alpha_0
      \ene
      Thus the above claim is true.  \\
      \indent For each $n$ with $-n<\alpha_0\leq\alpha_1<n$, $T_{n}$ is contained in $\bar{\Omega} \PLH [\alpha_0,\alpha_1]$ and $\bM(T_{n})$ is uniformly bounded and $\P T_n =\Gamma$. By Theorem \ref{compact:thm} as $n$ goes to $\infty$ the sequence $\{T_{n}\}_{n\geq \max\{|\alpha_0|,|\alpha_1|\}}$ will converge to a $T_\infty$ such that $ T_\infty|_{\bar{Q}_\phi}\subset\bar{ \Omega}\PLH [\alpha_0,\alpha_1]$ and 
      $$
      \bM(T_\infty):=\min \{\bM(T):\P T=\Gamma, T\in \mathcal{G}\}
      $$
      \indent Moreover from the compactness of BV functions in Theorem \ref{thm:compact:BV} there is a {\Ca} set $F_\infty$ such that $T_\infty=\P [[F_\infty]]|_{\bar{Q}_\phi}$ which is contained in $\bar{\Omega} \PLH[\alpha_0 ,\alpha_1 ]$ and $F_\infty$ is the subgraph of $\psi(x)$ outside $\bar{Q}_\phi$. By Theorem \ref{thm:decrease_perimeter}   there is a BV function $u(x)\in BV(\Omega')$ such that $u(x)=\psi(x)$ outside $\Omega$ and 
      \be 
       P_\phi (U , \Omega'\PLH \R)\leq  P_\phi (F_\infty, \Omega'\PLH \R)
       \ene 
       where $U$ is the subgraph of $u(x)$ and 
       \be 
       \alpha_0 \leq u(x) \leq \alpha_1 \text{ on } \Omega \quad u(x)=\psi(x)\text{  on  }  \Omega'\backslash \Omega
       \ene
        Moreover $\P(\P[[ U]]|_{\bar{Q}_\phi})=\Gamma$ by remark \ref{Rm:boundary}.  Since $F_\infty$ coincides with $U$ outside $\bar{Q}_\phi$, then 
       \be 
         \bM (\P [[U]] |_{\bar{Q}_{\phi}})\leq \bM(T_\infty)
       \ene 
       Thus $\P [[U]]|_{\bar{Q}_\phi}$ also solves the area minimizing problem \eqref{mass:minimizing}. This shows the existence of $u(x)$ with the property (2). \\
       \indent Now fix any $v(x)\in  BV(\Omega')$ satisfying $v(x)=\psi(x)$ outside $\bar{\Omega}$. Let $V$ be the subgraph of $v(x)$. By remark \ref{Rm:boundary} $\P(\P[[V]]|_{\bar{Q}_\phi})=\Gamma$. By the property (2)
       \be 
          \bM(\P [[U]] |_{\bar{Q}_{\phi}})\leq \bM(\P[[V]]|_{\bar{Q}_\phi})
       \ene 
       Because $u(x)=\psi(x)=v(x)$ outside $\bar{\Omega}$, by Theorem \ref{subgraph_and_area} and \eqref{eq:relation:st} we obtain 
       \be 
       \AF(u,\Omega')\leq \AF(v,\Omega')
       \ene 
      Thus $u(x)$ realizes the minimum of 
      \be 
         \min \{\AF(v(x),\Omega'): v(x)\in BV(\Omega'), v(x)=\psi(x)\text{  outside }\bar{\Omega} \}
      \ene 
       We conclude that $u(x)$ statisfies the property (1). \\
      \indent  At last we show that $u(x)$ has  the property (3). By the property (2) $U$ is a minimal set in $\Omega\PLH \R$. By Theorem \ref{thm:regularity:minimal:set},  except a closed set $S$ with $H^{n-6}(S)=0$, $\P U\backslash S$ is connected and $C^{1,\alpha}$ for some $\alpha\in (0,1)$. Since $\P U\backslash S$ is minimal in $\Omega\PLH\R\subset M_\phi$ and $\phi(x)$ is $C^2$, the regularity of minimal surface equation \eqref{eq:minimal:surface:equation:phi} implies that it is $C^{2,\beta}$ for any $\beta\in (0,1)$. Let $\vec{v}$ be the normal vector of $\P U\backslash S$ in $\Omega\PLH \R$ pointing to the positive infinity. Define $\Theta =\la \vec{v}, \P_r \ra$. Thus $\Theta\geq 0$ on $\P U\backslash S$ in $\Omega\PLH\R$ which is connected and $C^{2,\beta}$.  By the Harnack principle in corollary \ref{cor:harnack:principle}, we have $\Theta \equiv 0$ or $\Theta >0$ on $\P U\backslash S$. Since $\P U|_{\bar{Q}_\phi}$ is contained in $\bar{\Omega}\PLH[\alpha_0,\alpha_1]$ and $H^{n-1}(S)=0$,  only  $\Theta >0$ can happen in $\P  U\backslash S$. \\
      \indent Let $S'$ be the orthonormal projection of $S$ in $\Omega$. Therefore $H^{n-1}(S')=0$ and $u\in  C^{2,\beta}(\Omega\backslash S')$.  Moreover 
      \be 
      \int_{\Omega\backslash S'}\phi^{n}(x)\sqrt{1+|Du|^2}dvol \leq \AF(u,\Omega)<\infty 
      \ene 
      This implies that $u\in  W^{1,1}(\Omega)$. By the property (1) in  Theorem \ref{thm:area:minimizing} $u(x)$ is also the critical point of the functional 
      \be 
      \int_{\Omega}\phi^{n}(x)\sqrt{1+|Dv|^2}dvol\text{ for }v\in  W^{1,1}(\Omega)
      \ene 
      From the removable singularity result of Simon \cite[Theorem 1]{Sim77} $u(x)$ is regular at every $x$ in $\Omega$. Thus $u(x)$ satisfies \eqref{eq:minimal:surface:equation:phi} in the Lipschitz sense.  Since $\phi(x)\in C^2$, by the classical regularity of uniformly elliptic equations $u(x)$ is $C^{2,\beta}(\Omega)$ for any $\beta\in (0,1)$. The proof is complete. 
      \ep 
      
\section{the Dirichlet problem}
In this section we apply Theorem \ref{thm:area:minimizing} to solve the Dirichlet problem of minimal surface equations in $M_\phi$. 
\bd  We say that $\Omega$ is $\phi$-mean convex if the mean curvature of $\P\Omega$ satisfies $H_{\P\Omega}+n\la D \log\phi(x), \gamma\ra\geq 0$. Here  $D$ denotes the covariant derivative of $N$ and $H_{\P \Omega}=div(\gamma)$ for the outward normal vector $\gamma$ on $\P\Omega$. 
\ed
\br\label{Rm:phimeanconvex}  By \eqref{det:expression} (see also \cite[Lemma 3.1]{Zhou19-3}) $\Omega$ is $\phi$-mean convex if and only if $\Omega\PLH\R$ is mean convex in $M_\phi$ with respect to the metric $\phi^2(x)(\sigma+dr^2)$ for its outward normal vector. 
\er 
\subsection{The Dirichlet problem on $\phi$-mean convex domain} The first result of this section is given as follows. 
\bt\label{thm:DR:phimeanconvex} Let $\Omega$ be a $C^2$ bounded $\phi$-mean convex domain in $N$ and let $\phi(x)$ be a $C^2$ positive function in $N$. For any $\psi(x)\in C(\P\Omega)$ the Dirichlet problem 
\be \label{eq:DR:phimeanconvex}
\left\{
\begin{split}
L(u)=-div(\F{Du}{\omega})&+n\la  D \log \phi(x), -\F{Du}{\omega}\ra=0,  \text{ on } \Omega\\
        &u(x)=\psi(x)\quad \text{ on } \P\Omega
\end{split}\right.
\ene 
has a unique solution in $C(\bar{\Omega})\cap C^{2}(\Omega)$. Here $\omega=\sqrt{1+|Du|^2}$.
\et 
\br \label{Rm:CHHwork}
Suppose $\Omega$ is $C^{2,\alpha}$ for some $\alpha\in (0,1)$. Casteras-Heinonen-Holopainen \cite[Theorem 2]{CHH19} showed that if there is a positive constant $F>0$ such that \begin{enumerate}
	\item $m(x)\in C^2(N)$, $r(t)\in C^1(\R)$ satisfying 
	$$\max_{\bar{\Omega}}|Dm(x)|+\max_{t\in\R}|r'(t)|\leq F$$
	\item the Ricci curvature of $\Omega$ satisfies $Ric_\Omega\geq -\F{F^2}{n}$ and $H_{\P\Omega}\geq F$
	\end{enumerate} 
Then the Dirichlet problem 
\be 
-div(\F{Du}{\omega})+\la D(m(x)), -\F{Du}{\omega}\ra +r'(u(x))\F{1}{\omega}=0\text{ on } \Omega 
\ene 
with $u(x)=\psi(x)$ on $\P\Omega$ for any $\psi(x)\in C(\P\Omega)$ has a solution in $C^{2,\alpha}(\Omega)\cap C(\bar{\Omega})$. \\
\indent 
 In the case of $r(t)\equiv 0$ Theorem \ref{thm:DR:phimeanconvex} removes the curvature assumption on $\Omega$ in \cite[Theorem 2]{CHH19}. Our mean curvature assumption on $\P\Omega$ should be optimal. For example when $\phi(x)\equiv 1$ this is confirmed by Serrin \cite{Serrin69} in Euclidean spaces. 
\er 
\br\label{rk:continuous:method}  The above result can not be obtained by the continuous method in \cite[section 11.2, 11.3,  Chapter 18]{GT01}. The reason is that the boundary assumption may not be preserved by these methods. 
\er 
\bp The uniqueness of the Dirichlet problem \eqref{eq:DR:phimeanconvex} is obvious by the maximum principle. We only need to show the existence.  \\
\indent First we assume that $\psi(x)\in C^3(\P\Omega)$.  Without loss generality we can assume $\psi(x)\in C^1(\Omega'\backslash \Omega)$ for some $\Omega'$ strictly containing $\Omega$. Let $Q_\phi$ be the set $\Omega\PLH \R$ and $\bar{Q}_\phi$ be its closure. By Theorem \ref{thm:area:minimizing} there is a $u(x)\in BV(\Omega')$ with $u(x)=\psi(x)$ outside $\Omega$ such that $T=\P[[U]]|_{\bar{Q}_\phi}$ realizes the area minimizing problem \eqref{mass:minimizing}. Moreover $u(x)\in C^2(\Omega)$. As a result $\P U\cap Q_\phi$ is a minimal graph in $Q_\phi$. By the conclusion (1) in Theorem \ref{thm:minimal:equation} $u(x)$ satisfies $Lu=0$ on $\Omega$. It only suffices to show that $u(x)\in C(\bar{\Omega})$\\
\indent Fix any $z\in \P\Omega$. Now define 
\be 
A:=sup\{\lim_{n} \sup_{k\geq n} u(x_k): \{x_n\}_{n=1}^\infty\in \Omega, \lim_{n\rightarrow \infty}dist(x_n, z)=0 \}
\ene 
where $dist$ is the distance function in $N$.  By Theorem \ref{thm:area:minimizing} $ \min_{\P\Omega}\psi \leq u(x)\leq \max_{\P\Omega}\psi$ for any $x\in \Omega$. Thus $A$ is a finite number. Suppose $A>\psi(z)$. There is a sequence $\{x_n \}_{n=1}^\infty$ in $\Omega$ such that $\lim_{n\rightarrow +\infty} x_n=z$ and $\lim_{n\rightarrow +\infty}\psi(x_n)=A$. Since $\psi(x)$ is continuous at $z$, there is a neighborhood $V$ of the point $(z, A)$ such that $V$ is disjoint with $\Gamma$, the graph of $\psi(x)$ in $\P\Omega\PLH\R$.\\
\indent By Theorem \ref{thm:area:minimizing} $U$ is a minimal set in $V\cap \bar{Q}_\phi$ and passes through the point $(z,A)$. Since $\P \Omega\PLH \R$ is $C^2$, Theorem \ref{regularity:key:thm} implies that $\P U$ is a $C^{1,\alpha}$ graph near $(z,A)$.  Since $U$ is a minimal set in $V\cap \bar{Q}_\phi$, then $H_{\P U}\leq 0$ near $(z,A)$ with respect to the outward normal vector of $\Omega\PLH\R$ in the Lipschitz sense. Since $\Omega$ is $\phi$-mean convex, by remark \ref{Rm:phimeanconvex} $\Omega\PLH \R$ is a mean convex domain in $M_\phi$. Because $ U\cap V\subset \Omega\PLH \R$ and  is tangent to $\P\Omega\PLH\R$. From the maximum principle in  \cite[Theorem A.1]{GZ20}, $\P U$ coincides with $\P\Omega \PLH \R$ in $V$. This contradicts to the definition of $A$. Thus 
\be\label{def:argument:one}
\lim \sup u(x_n)\leq \psi(z), \forall \{x_n\}_{n=1}^\infty\in \Omega\,\lim_{n\rightarrow \infty }x_n=z 
\ene 
With a similar derivation we also obtain 
\be\label{def:argument:two}
\lim \inf u(x_n)\geq \psi(z), \forall \{x_n\}_{n=1}^\infty\in \Omega\, \lim_{n\rightarrow \infty }x_n=z
\ene 
Combining the above two facts together yields that $u(x)$ is continous for any fixed $z\in \P\Omega$, thus $u(x)\in C(\bar{\Omega})$. As a result we obtain the existence of the Dirichlet problem \eqref{eq:DR:phimeanconvex} when $\psi(x)\in C^3(\P\Omega)$.\\
\indent Now we consider the general case of $\psi(x)\in C(\P\Omega)$.  Then there are two sequences $\{\psi^{\pm}_k(x)\}_{k=1}^\infty$ in $C^3(\P\Omega)$ such that 
\be 
\cdots \psi^-_{k}(x)\leq \psi^-_{k+1}(x)\leq \cdots \psi(x)\leq \cdots \leq \psi^+_{l+1}(x)\leq \psi^+_{l}(x)\leq \cdots
\ene
and both of them converge to $\psi(x)$ in the $C^0(\P\Omega)$.  By the previous argument for any positive integer $k$ let $u^{\pm}_k(x)$ be the solution of the Dirichlet problem \eqref{eq:DR:phimeanconvex} on $\Omega$ with $u_k(x)=\psi^\pm_k(x)$ on $\P\Omega$. Applying the maximum principle on \eqref{eq:DR:phimeanconvex}, we obtain
\begin{enumerate} 
	\item $\max_{x\in \bar{\Omega}}|u_k(x)|\leq C$ where $C=\max_{k=1,\cdots,x\in\P \Omega} \{\psi(x),\psi^\pm_{k}(x)\}$;
	\item  $ u^-_{k}(x)\leq u^-_{k+1}(x)\leq \cdots \leq u^+_{l+1}(x)\leq u^+_{l}(x)\cdots $ on $\Omega$;
\end{enumerate}
According to Theorem \ref{thm:interior:estimate}, $\{u^\pm_k(x)\}_{k=1}^\infty$ has a local uniformly bound of their $C^2$ norms in $\Omega$. Now we can choose  a subsequence, still written as $\{u^-_k(x)\}_{k=1}^\infty$ such that $\lim_{k\rightarrow\infty} u^-_k(x)=u(x)$ on $\bar{\Omega}$ and $u(x)\in C^2(\Omega)$ and satisfies $Lu=0$ on $\Omega$. \\
\indent The only thing to left is to show that $u(x)\in C(\bar{\Omega})$ and $u(x)=\psi(x)$ on $\P\Omega$. By the conclusion (2) above applying the maximum principle on $Lu=0$ yields that 
\be 
   u^-_k(x)\leq u(x)\leq u^+_l(x)\text{ on } \Omega
\ene 
for any positive integers $k,l$.  This implies that for any $z\in \P\Omega$
\be 
\psi^-_k(z)\leq \lim_{x\rightarrow z,x \in \Omega}\inf u(x)\leq\lim_{x\rightarrow z,x \in \Omega}\sup u(x)\leq \psi^+_l(z)
\ene 
Letting $k,l$ go to the positive infinity we obtain $\lim_{x\rightarrow z, x\in \Omega} u(x)=\psi(z)$ for any $z\in\P\Omega$.  Thus we define $u(z)=\psi(z)$ on $\P\Omega$ and $u(x)\in C(\bar{\Omega})$. The proof is complete. 
\ep 
\subsection{The uniqueness of area minimizing currents}
Now by Theorem \ref{thm:DR:phimeanconvex} we shall obtain a uniqueness result for the area minimizing problem \eqref{mass:minimizing}. 
\bt\label{thm:uniqueness:mass:minimizng} Suppose $\Omega$ is a $\phi$-mean convex $C^2$ domain and $\psi(x)\in C^1(\P\Omega)$. Set $\Gamma =\{(x,\psi(x)):x\in\P\Omega\}$ . Let $u(x)$ be the solution of the Dirichlet problem \eqref{eq:DR:phimeanconvex} with boundary data $\psi(x)$.  Let $U(x)$ be the subgraph of $u(x)$. 
\begin{enumerate}
	\item  Then $T=\P [[U]]|_{\bar{Q}_\phi}$ is the unique n-{\IMC} in $\mG$ to realize 
\be \tag{\ref{mass:minimizing}}
\min\{\bM(T): T\in  \mG, \P T=\Gamma \}
\ene 
\item Let $\Omega'$ be a domain with $\Omega\subset\subset \Omega'$. Extend $\psi(x)$ as a function $C^1(\Omega'\backslash \Omega)$. Let $u(x)=\psi(x)$ outside $\Omega$. Then $u(x)$ is a unique function to realize 
\be \tag{\ref{functional:area:minimizing}}
\min \{\AF(v(x),\Omega'): v(x)=\psi(x) \text{  outside $\Omega$}\}
\ene
\end{enumerate}
\et 
\br  This theorem is a finite version of \cite[Theorem 4.1]{HL87} in finite conformal cones in $M_\phi$. 
\er 
\bp By Theorem \ref{thm:DR:phimeanconvex}, let $T$ be a n-{\IMC} with compact support in $\mathcal{G}$ to realize \eqref{mass:minimizing}. For any $t\in \R$, let $U_t$ be the subgraph of $u(x)+t$ over $\Omega$. By \eqref{eq:DR:phimeanconvex}, $\P[[ U_t]]$ is a minimal graph over $\Omega$ in $M_\phi$  and the boundary of $\P [[U_t]]$ is disjoint with $\Gamma$ if $t\neq 0$. \\
\indent Now the following number
\be 
t_0=\inf\{t>0:  \P [[U_t]]\cap T =\emptyset:  \forall \quad s \in [t, +\infty)\}
\ene 
is well-defined since $spt(T)$ is compact. Suppose $t_0>0$. Then $\P [[U_{t_0}]]$ is tangent to $T$ at some point $p\in \Omega\PLH\R$. \\
\indent By Lemma \ref{lm:au:mass:minimizing}, $T=(\P [[F]])|_{\bar{Q}_\phi}$ where  $F$ is a minimal set in the set $\bar{\Omega}\PLH\R$. From the definition of $t_0$, $F\cap \Omega\PLH \R \subset U_{t_0}$. By Theorem \ref{regularity:key:thm} $\P F$ is $C^{1,\alpha}$ near $p$. In the following let $H$ denote the mean curvature. Thus $H_{\P F}\leq 0$ near $p$ with respect to the upward normal vector in the Lipschitz sense. Note that  $H_{\P U_{t_0}}=0$ near $p$.  Then the maximum principle in  \cite[Theorem A.1]{GZ20} implies that $\P F$ coincides with $\P U_{t_0}$ near $p$. By the connectedness of $\P U_{t_0}$,  we have $\P U_{t_0}\subset  T$. Thus for  some $x_0\in \P\Omega$, the point $p_0=(x_0,u(x_0)+t_0)$ is contained in $T$. \\
\indent Since $t_0>0$,  there is a neighborhood $V$ of $p_0$ disjoint with $\Gamma$. Then $F$ is a minimal set in $V\cap \bar{\Omega}\PLH \R$.  By Remark \ref{Rm:phimeanconvex}, $\Omega\PLH\R$ is mean convex in $M_\phi$. Again by \cite[Theorem A.1]{GZ20}, $\P\Omega \PLH\R$ coincides with $T$ near $p_0$. With the connectedness, the set $\{(x,r): x\in\P\Omega, r>u(x)+t_0\}$ is contained in $T$. This is impossible since $spt(T)$ is compact. Thus  $t_0=0$ and $spt(T)\subset \{(x,r):x\in \bar{\Omega}, r\leq u(x)\}$.  With a similar derivation, we can show that $spt(T)\subset \{(x,r):x\in \bar{\Omega}, r\geq u(x)\}$. Finally combining the above two facts yields that $T=\P  [[U]]|_{\bar{Q}_\phi}$ where
 $U$ is the subgraph of $u(x)$. We obtain the conclusion (1). \\
\indent Suppose there are two functions $u_1, u_2$ in $BV(\Omega')$ to solve \eqref{functional:area:minimizing}. Let $U_i$ be the subgraph of $u_i(x)$ and $T_i=\P [[U_i]]|_{\bar{Q}_\phi}$ for $i=1,2$. By (2) in Theorem \ref{thm:area:minimizing}, $T_i$ should solve \eqref{mass:minimizing}. By the conclusion (1) we just obtained, $T_1=T_2$ and $u_1(x)=u_2(x)$ on $\Omega$.  We arrive the conclusion (2). The proof is complete. 
\ep 
\section{Infinity boundary cases} 
In this section we apply the results in previous sections to generalize the existence and the uniqueness of area minimizing graphs with the infinity star-shaped boundary in hyperbolic spaces. \\
\indent Throughout this section, we assume that $N$ is a compact $n$-dimensional Riemannian manifold with its metric $\sigma$ and compact embedded $C^2$ boundary $\P N$. We define 
\be \label{assumption:one}
N_r:= \{d(x)>r: x\in N\}
\ene 
where $d(x)$ denotes $d(x,\P N)$, the distance between $x$ and $\P\Omega$. The following lemma is obvious. 
\bl \label{def:ro}
 Since $\P N$ is $C^2$ embedded and compact, there is a $ r_0>0$ such that for any $x\in N\backslash N_{r_0}$, there is a unique $y \in \P N$ such that  $d(x)$ is equal to $d(x,y)$.  
 \el
 From now on fix $r_0$. Thus $d(x)$ is $C^2$ on $N\backslash N_{r_0}$. Suppose $\phi(x)$ is a positive $C^3$ function on $N$ with 
\be\label{assumption:two}
\phi(x)=h(d(x)) \text{  on } N\backslash N_{r_0}\quad  \lim_{d(x)\rightarrow 0} \phi(x)=\lim_{d(x)\rightarrow 0} h(d(x))=+\infty 
\ene  
where $h(r):(0,r_0)\rightarrow \R^+$ is a positive $C^2$ function. Then 
\be M_\phi:=(N\PLH, \phi^2(\sigma+dr^2)))
\ene is a complete Riemannian manifold. A natural compactification is 
\be 
\bar{M}_\phi =M_\phi\cup (\P N\PLH\R)
\ene 
equipped with the product metric topology of $(N\PLH\R, \sigma+dr^2)$. 
\bd Let $S$ be a complete $k$-{\IMC} in $M_\phi$.  Its infinity asymptotic boundary in $\P N\PLH\R$ is the set $\bar{S}\backslash S$ where $\bar{S}$ is the closure of $S$ in the product metric topology. 
\ed 
\bd Let $\bM$ be the mass of currents in $M_\phi$. We say that $T$ is a local area minimizing $n$-{\IMC} in $M_\phi$ if $spt(T)$ is contained in $N\PLH [a,b]$ for some finite interval $[a,b]$, for any $n$-{\IMC} $T'$ satisfying $T=T'$ outside a compact set in $M_\phi$, then $\bM(T)\leq \bM(T')$. 
\ed 
\indent  The main result of this section is given as follows. 
\bt \label{thm:generalization:Hyperbolic} Let $N_r,\phi, r_0$ be given in \eqref{assumption:one}, \eqref{assumption:two} and Lemma \ref{def:ro} respectively. Suppose for any $r\in (0, r_0)$ $N_r$ is $\phi$-mean convex. Fix any $\psi(x)\in C^1(\P N)$.  Let $\Gamma$ denote $\{(x,\psi(x)):x\in\P N\}$ in $\P N\PLH\R$. Then there is a unique local area minimizing $n$-{\IMC} $T$ in $M_\phi$ with the infinity boundary $\Gamma$. Moreover $T$ is a graph over $N$.
\et  
\br\label{Rm:hyperbolic:space} Recall that 
$\mathbb{H}^{n+1}$ is the upper half space $\{(x,y):x\in \R^n, y>0\}$ equipped with the metric
\be 
ds^2_H= \F{ dx^2+dy^2}{y^2}\quad \text{ $dx^2$ is the Euclidean metric in $\R^n$}
\ene 
Let $S^n_{+}$ be the upper half hemisphere in $\R^{n+1}$ with the induced metric $\sigma_n$. We introduce the pole coordinate on $S^n_+$ (the pole at the north pole) $(\theta,\vp)\in S^{n-1}\PLH [0,\F{\pi}{2}]$. Thus $\sigma_n=d\vp^2+\sin^2(\vp) d\theta^2$. Note that $dx^2+dy^2=s^2 \sigma_n+ds^2$ for $s\in (0,\infty)$ and $x_{n+1}=s \cos(\vp)$. Letting $r=\ln s$ we can represent Hyperbolic spaces $\mathbb{H}^{n+1}$ as 
\be 
(S^n_+\PLH \R,  \F{1}{\cos^2(\vp)}(\sigma_n+dr^2))
\ene 
which is a special case of \eqref{eq:Model}. \\
\indent The above theorem is to generalize Hardt-Lin's result in Hyperbolic spaces $\mathbb{H}^{n+1}$ \cite[Theorem 4.1]{HL87}. We verify this claim as follows. Fixing any number $\vp_0\in (0, \F{\pi}{2})$ we define a domain 
 \be 
 S_{\vp_0}:=\{(\theta,\vp)\in S^n_+: \vp\in [0, \vp_0)\}
 \ene
 
By \cite[Proposition 2.1]{Zhou18} the mean curvature of $\P S_{\vp_0}$ {\wrt} $\F{\P}{\P \vp }$ is $(n-1)\F{\cos(\vp_0 )}{\sin(\vp_0)}$. By \eqref{det:expression} the mean curvature of $\P S_{\vp_0}\PLH \R$ in $\mathbb{H}^{n+1}$ is 
\be\label{S_k:mean:curvature}
\begin{split}
H_{ \P S_{\vp_0}\PLH \R}&=\cos(\vp)\{(n-1)\F{\cos(\vp)}{\sin(\vp)}+n\la D\log\F{1}{\cos(\vp)}, \F{\P}{\P \vp}\ra \}|_{\vp=\vp_0}\\
                 &=\cos(\vp_0)\{(n-1)\F{\cos(\vp_0)}{\sin(\vp_0)}+n\F{\sin(\vp_0)}{\cos(\vp_0)} \}
\end{split}
\ene 
Thus $S_{\vp_0}\PLH \R$ is mean convex in $M_\phi$ and $S_{\vp_0}$ is $\phi$-mean convex by Remark \ref{Rm:phimeanconvex}.  Then $\mathbb{H}^{n+1}$ satisfies the condition in the above theorem. Thus  \cite[Theorem 4.1]{HL87} is a speical case of Theorem \ref{thm:generalization:Hyperbolic}.  
 \er 
 Now we are ready to show Theorem \ref{thm:generalization:Hyperbolic}. 
 \bp  First we show the existence of a minimal graph in $M_\phi$ with the infinity boundary $\Gamma$. The local area minimizing property can be easily obtained by Theorem \ref{thm:area:minimizing}.  \\
 \indent \textit{Step one: Existence.} First assume $\psi(x)\in C^3(\P N)$. Since $\P N$ is compact, by Lemma \ref{def:ro} we extend $\psi(x)$ into a $C^2$ function on $N$ such that 
 \be \label{eq:extension:psi}
 \psi(x)=\psi(y)\text{ for any $x\in N\backslash N_r$ with $d(x,y)=d(x,\P N)$}
 \ene 
   Because $N_r$ for each $r<r_0$ is $\phi$-mean convex,  by Theorem \ref{thm:DR:phimeanconvex} there is a $u_r(x)$ in $C^2(\Omega)\cap C(\bar{\Omega})$ to solve the Dirichlet problem 
 \be \label{mid:eq:DR:phimeanconvex}
 \left\{
 \begin{split}
 Lu:=-div(\F{Du}{\omega})&+n\la D\log \phi(x), -\F{Du}{\omega}\ra=0,  \text{ on } N_r\\
 	&u(x)=\psi(x)\quad \text{ on } \P N_r
 \end{split}\right.
 \ene 
 Let $\mu$ be a positive constant such that $\max_{x\in N\cup \P N}|\psi(x)|\leq \mu$. By the maximum principle, for each $r\in (0,r_0)$, 
 \be\label{def:mu}
 \max_{N_r}|u_r|\leq \mu
 \ene 
 For any fixed embedded open ball $B_{\Sc}(x)$ with $dist(x,\P N)>2\Sc$, Theorem \ref{thm:interior:estimate} implies that 
 \be \label{def:mu1}
\max_{\bar{B}_\Sc(x)}\{|u_r|(x), |Du_r|(x), |D^2u_r|(x)\}\leq \mu_1
 \ene 
 where $\mu_1$ is a constant depending on $\psi$, $\phi$ on $B_{\F{3}{2}\Sc}(x)$ and $\mu$. Thus after choosing a subsequence from $\{u_r(x)\}_{r<r_0}$, denoted by $u_k(x)$, $u_k(x)$ converges uniformly to $u_\infty(x)$ in the $C^2$ sense on any compact set of $N$ disjoint with $\P N$. Thus 
 \be
  Lu_\infty(x)=0 \quad \text{  on  }  N
 \ene 
 Next we show $\lim_{x\rightarrow z,x \in N}u_\infty(x)=\psi(z)$ for any $z\in \P\Omega$. To achieve this goal we construct the supersolution and subsolution of $Lu=0$ in \eqref{mid:eq:DR:phimeanconvex} on $N_{r}$ for some $r\in (0, r_0)$. 
 \bl\label{lm:supper:solution}  Let $\psi(x)$ satisfy \eqref{eq:extension:psi}, let $\mu$ be a constant given in \eqref{def:mu} and $r_0$ be given in Definition \ref{def:ro}. Suppose for any $r\in (0, r_0)$ $N_r$ is $\phi$-mean convex. Then there are three positive constants $r_1\in (0, r_0)$ and $\nu>0$ and $\kappa>0$ such that 
 \begin{gather}
u_{\pm}(x)=\psi(x)\pm \vp(d(x))  \text{ on $N\backslash N_{r_1}$\quad }\vp(r):=\F{1}{\nu}\log(1+\kappa r)\label{eq:def:vp} \\
L u_+(x)\geq 0, \quad L u_{-}(x) \leq 0 \text{ on $N\backslash N_{r_1}$}\label{conclusion:one}\\
u_+(x)\geq 2\mu, \quad u_{-}(x)\leq -2\mu \text{ on $\{x\in N:d(x,\P N)=r_1\}$}\label{conclusion:two}
 \end{gather}
\el  
\bp Let $\{x_1,\cdots, x_n\}$ be a local coordinate in $N$ and $\{\P_i\}_{i=1}^\infty$ be the corresponding frame. We write $\la \P_i,\P_j\ra$ as $\sigma_{ij}$ and $(\sigma^{ij})=(\sigma_{ij})^{-1}$.  For any smooth function $f$,  $f_i$, $f_{ij}$ denote the first and second covariant derivatives of $f$ and $f^i=\sigma^{ik}f_k$ with the sum over $k$.\\
\indent Note that \be 
\mathcal{E} u:=\omega^3 Lu=(1+|Du|^2)\{-\Delta u+ n\la D \log\phi(x), -Du\ra\}+u^i u^j u_{ij}
\ene 
Here $\omega=\sqrt{1+|Du|^2}$,  $\Delta u=div(Du)$ and $D$ denote the Laplacian and the covariant derivative of $N$ respectively. Now assume $u_{\pm}(x)$ are given in \eqref{eq:def:vp} where $\nu, \kappa$ are determined later. \\
\indent By \eqref{assumption:two} and \eqref{eq:extension:psi} on $N\backslash N_{r_0}$ we have $\la D\psi, Dd\ra = 0$ and 
\be \label{eq:fact:two}
\la D\log \phi(x), -Du_+(x)\ra=\vp'\la D\log \phi(x), -Dd\ra\quad  
\ene  where $d(x)=d(x,\P N)$ and $\vp'=\F{1}{\nu(1+\kappa d(x))}$. With \eqref{eq:fact:two} we compute $\mathcal{E}u_+$ on $N\backslash N_{r_0}$ as follows. 
\begin{align*}
\mathcal{E}(u_+) &=\vp' (1+|Du_+|^2)(-\Delta d+n\la D\log\phi, -Dd\ra )-(1+|D\psi|^2)\Delta \psi\\
&+\psi^i\psi^j\psi_{ij}
+(\vp')^2(-\Delta \psi+d^i d^j\psi_{ij})-(1+|D\psi|^2)\vp''
\end{align*}
 Since $N_r$ is $\phi$-mean convex for each $r\in (0, r_0)$ and $\vp'>0$, the first term above is nonnegative. By \eqref{eq:def:vp} $\vp''=-\nu\vp'^2$. Now assuming $\vp'\geq 1$ $\mathcal{E}(u_+)$ on $N\backslash N_{r_0}$ satisfies 
\be 
\mathcal{E}u_+\geq \nu \vp'^2-C(\vp'^2)
\ene 
for some positive constant $C$ only depending on $\psi(x)$. Now take $\nu= C$ we obtain $\mathcal{E}u_+\geq 0$ on $N\backslash N_{r_0}$. This shows \eqref{conclusion:one} for $u_+$ and we have to determine $\kappa$ and $r_1\in (0, r_0)$ in a such way $\vp'(r)\geq 1$ and $\vp(r_1)\geq \mu_1:= 2\mu+ \max_{\bar{N}}\psi(x)$. We have 
\begin{gather}
\vp'(r)=\F{1}{\nu}\F{\kappa}{1+\kappa r}>\F{1}{\nu}\F{\kappa}{1+\kappa r_1}\geq 1\\
\vp(r_1)=\F{1}{\nu}\log (1+\kappa r_1)\geq \mu_1
\end{gather}
All conditions are satisfied provided we take $r_1$ small enough and large $\kappa$ satisfying \be 
\kappa \geq \max\{ \F{\nu}{1-r_1\nu}, \F{e^{\mu_1\nu}}{r_1}\}
 \ene  This gives \eqref{conclusion:two} for $u_+(x)$. \\
\indent A similar derivation yields the conclusions for $u_-(x)$. The proof is complete. 
\ep 
Now we continue to show Theorem \ref{thm:generalization:Hyperbolic}. 
Now for any $r\in (0,r_1)$, with the maximum principle, \eqref{mid:eq:DR:phimeanconvex} and \eqref{def:mu} Lemma \ref{lm:supper:solution} implies that $u_-(x)\leq u_r(x)\leq u_+(x)$ on $N_r \backslash N_{r_0}$. Let $r$ go to $0$ we obtain that 
\be 
 u_-(x)\leq u_\infty(x)\leq u_+(x)\quad \text{ on $N \backslash N_{r_1}$}. 
\ene 
Since both $u_+(x)$ and $u_-(x)$ are continuous on $N\cup \P N$ and equal to $\psi(x)$ on $\P N$, thus 
\be \label{eq:boundary:continouity}
\lim_{x\in N, x\rightarrow z\in \P N} u_\infty(x)=\psi(z)
\ene  Let $U$ be the subgraph of $u_\infty(x)$ in $M_\phi$. Let $T=\P [[U]]$ be the corresponding {\IMC}. Thus with respect to the product topology of $N\PLH\R$, $\P T=\Gamma$. \\
\indent Now we obtain the existence of a local {\IMC} in $M_\phi$ with the desirable infinity boundary $\Gamma$ for any $\psi \in C^3(\P N)$. Moreover it is a minimal graph over $N$ in $M_\phi$.\\
\indent \textit{Step Three: General Case} As for $\psi(x)\in C(\P N)$, we can construct  two monotone sequences $\{ \psi^1_k(x)\}_{k=1}^\infty$ and $\{ \psi^2_k(x)\}_{k=1}^\infty$ in $C^3(\P N)$ such that the former one converges increasingly to $\psi(x)$ and the latter one converges decreasingly to $\psi(x)$ in the $C^0(\P N)$ sense. Let $\{u^i_{k}(x)\}_{k=1}^\infty$ be the solutions of $L u=0$ on $N$  with asymptotic value $\{\psi^i_k(x)\}_{k=1}^\infty$ for $i=1,2$. From the maximum principle and the translating invariant of minimal graphs, we have 
$\{u^1_k(x)\}$ is an increasing sequence on $N$ and $u^1_k(x)\leq u^2_l(x)$ for any $k,l$. Note that both sequences are uniformly bounded. By the interior estimate of $Lu=0$ in $N$, $u_k(x)$ locally converges to a $C^2$ function $u_\infty(x)$ in $N$ satisfying $Lu_\infty=0$.  
Thus for any $k,l$, on $N$ we obtain 
\be 
   u^1_k(x)\leq u_\infty(x)\leq u^2_l(x)
\ene 
Thus \eqref{eq:boundary:continouity} still holds for $\psi(x)\in C(\P N)$.  As a result for any $\psi(x)\in C(\P N)$ we show there is a minimal graph $\Sigma=(x,u(x))$ in $M_\phi$ with the infinity boundary $\Gamma$. It is also a local area minimizing {\IMC} according to Theorem \ref{thm:area:minimizing}. \\
\indent \textit{Step Four: Uniqueness.}  Suppose $T$ is a local area minimizing $n$-{\IMC} with infinity boundary $\Gamma=(x,\psi(x))$. For any $t\in \R$, we define $f_t(x,t)=(x,r+t)$ for any $x\in N$ and $t\in \R$. By \eqref{eq:Model}, $f_t$ is an isometry of $M_\phi$. Thus $f_{t,\#}T$ is also a local area minimizing $n$-{\IMC} in $M_\phi$ with infinity boundary $\Gamma_t=(x,\psi(x)+t)$. Recall that $spt(T)\subset N\PLH [a,b]$. We claim that $f_{t,\#}(T)\cap T=\emptyset$ for any $t\neq 0$. Otherwise there is a $t\neq 0$ such that the regular part of $f_{t,\#}(T)$ and $T$ intersects transversely. Since $\Gamma_t\cap \Gamma =\emptyset$ for any $t\neq 0$. By the area minimzing property, we can replace a piece of $T$ by the piece of $f_{t,\#}$ that is cut off by $T$ and still have a local area minimizing property. But then the singular set of the resulting {\IMC} would include the intersection of $T$ and $f_{t,\#}(T)$ with $n-1$ dimension. This contradicts to the $n-7$ dimensional singular set of local area minimizing {\IMC}s. Thus $T$ is a graph over $N$. Arguing similar as in the proof of the property (3) in Theorem \ref{thm:area:minimizing}, $T$ is a $C^2$ minimal graph in $M_\phi$ with the infinity boundary $\Gamma$. Since $N\cap \P N$ is compact, applying the maximum principle into $Lu=0$ yields that the uniqueness of the minimal graph. \\
\indent  The proof is complete. 
 \ep 
\appendix 
\section{Interior estimate of mean curvature equations} \label{interior-estimate-}
Throughout this section let $N$ be a complete Riemannian manifold with a metric $\sigma$ and $f(x)$ be a $C^3$ function in $N$. Let $M_0$ be the product manifold $N\PLH\R$ equipped with the metric $\sigma+dr^2$. Let $D$ be the covariant derivative of $N$.  Let $S$ be an orientable $C^{2}$ hypersurface in $M_0$ with the normal vector $\vec{v}$ for any $\beta_1\in (0,1)$.  We have the following theorem about its angle function $\la \vec{v}, \P_r\ra $. Comparing to  \cite[section 2.3]{CHH19} we use a moving frame method from \cite[section 2]{Eich09} to obtain the interior estimate of mean curvature equations. 
\bl\label{lm:angle} Let $\Theta=\la \vec{v},\P_r\ra$. Suppose the mean curvature of $S$ with respect $\vec{v}$ satisfies that 
\be \label{eq:mean:curvature:assumption}
H_S+\la Df, \vec{v}\ra=0
\ene 
where $D$ is the gradient of $f$ on $N$. Then in the Lipschitz sense it holds that   
\be\label{formula:theta}
\begin{split}
\Delta \Theta &+(|A|^2+\bar{R}ic(\vec{v},\vec{v}) -Hess(f)(\vec{v},\vec{v}\ra)\Theta \\
&+ \la \nb \Theta, Df\ra =0      
\end{split}
\ene 
where $|A|^2$ is the second fundamental form of $S$, $\nb$ is the covariant derivative of $S$ and $\bar{R}ic$ is the Ricci curvature of the product manifold $M_0$, $Hess$ is the Hessian of $f$ in $M_0$. 
\el 
\bp  
Notice that $S$ is a $C^3$ hypersurface. We denote $\{e_i\}_{i=1}^n$ by a local orthonormal frame of $S$. Note that the derivation in Lemma 2.2, \cite{Zhou19} is true for any dimension. That is by Lemma 2.2 in \cite{Zhou19} we have 
\be \label{eq:first:evolution}
\Delta \Theta+(|A|^2+\bar{R}ic(\vec{v},\vec{v}))\Theta-\la \nb H, \P_r\ra=0
\ene 
where $H$ is the mean curvature of $S$. Let $\bar{D}$ is the corvariant derivative of $M_0$. Since $-H= \la Df, \vec{v}\ra$, 
\be \label{de:mid:one}
-\la \nb H, \P_r\ra =\la e_i,\P_r\ra \la \bar{D}_{e_i}Df, \vec{v}\ra+ \la e_i,\P_r\ra \la Df, e_k\ra h_{ik}
\ene 
where $\bar{D}_{e_i}\vec{v}=h_{ik}e_k$. Note that $\la \P_r, e_i\ra e_i=\P_r-\la \vec{v}, \P_r\ra \vec{v}$. Thus 
\be \label{de:mid:two}
\la e_i,\P_r\ra \la \bar{D}_{e_i}Df, \vec{v}\ra = -\Theta Hess(f)\la \vec{v},\vec{v}\ra
\ene 
Here we use the fact $\bar{D}_{\P_r }Df =0$ . Since the metric of $M_0$ is the product metric, $\bar{D}_X \P_r =0$ for any tangent vector field $X$.  Thus $\la \nb \Theta, Df\ra= \la Df, e_k\ra h_{ik}\la e_i,\P_r\ra$. Combining this with \eqref{de:mid:one}, \eqref{de:mid:two} we obtain 
\be 
-\la \nb H, \P_r\ra =-\Theta Hess(f)\la \vec{v},\vec{v}\ra+\la \nb \Theta, Df\ra
\ene 
Putting this into \eqref{eq:first:evolution} we obtain the conclusion in the case $f\in C^3(N)$. \\
\indent In the case of $f\in C^2(N)$, $S$ is a $C^{2,\beta}$ hypersurface for any $\beta$ in $(0,1)$.  The result comes from the classical approximating method.  The proof is complete. 
 \ep 
 As a corollary we obtain the following Harnack type result. 
 \begin{cor} \label{cor:harnack:principle} Let $\Omega$ be a bounded domain in $N$ and $f$ be a $C^3$ function in $N$. Suppose $S$ is a $C^2$ connected orientable hypersurface satisfying \eqref{eq:mean:curvature:assumption} with its normal vector $\vec{v}$.  If $\Theta=\la \vec{v},\P_r\ra \geq 0$ on $S$, $\Theta\equiv 0$ or $\Theta >0$ on the whole $S$. 
 \end{cor}
\bp Suppose $f$ is $C^3$. Then $S$ is $C^3$. Because $\Omega$ is bounded,  there is a positive constant $C>0$ such that 
\be \label{harnack:inequality}
\Delta \Theta -C\Theta + \la \nb \Theta, Df\ra \leq 0     
\ene 
Since $\Theta\geq 0$ on $S$, the weak maximum principle implies that $\Theta\equiv 0$ or $\Theta >0$ on the whole $S$. \\
\indent When $f$ is $C^2$, $S$ is $C^{2,\beta}$.  Then $\Theta$ satisfies \eqref{harnack:inequality} in the Lipschitz sense. By the weak Harnack inequality in \cite{Tru67}, $\Theta>0$ or $\Theta\equiv 0$ on the whole $S$. The proof is complete. 
\ep 
We have the following interior estimate of mean curvature equations. 
\bt \label{thm:interior:estimate}
Let $f\in C^3(N)$. Fix $x_0\in N$ and let $B_\rho(x_0)$ be an embedded ball centered at $x_0$ with radius $\rho$.  Suppose $u(x)\in C^2(B_\rho(x_0))$ satisfying 
\be \label{eq:mc:def}
-div(\F{Du}{\omega})+\la Df, -\F{Du}{\omega}\ra=0
\ene 
where $Du$ is the gradient of $u$, $\omega=\sqrt{1+|Du|^2}$ and $div$ is the divergence of $N$. 
If $ \max_{x\in B_\rho(x_0)}|u(x)|\leq c_0$ for some positive constant $c_0$, then 
\be 
\max_{B_\F{\rho}{2}(x_0)}|Du|\leq C
\ene 
where $C$ is a constant only depending on the Ricci curvature, the $C^2$ norm of $f$, $c_0$ and $\rho$. 
\et 
\bp Our proof follows from the idea of Lemma 2.1 in \cite{Eich09} by Eichmair. \\
\indent  First we assume $f(x)$ is a $C^3$ function in $N$. By the classical Schauder estimate $(x)$ is $C^{3,\beta}$ for any $\beta\in (0,1)$.  We denote by $\Lambda$ the $C^2$ norm of $f$ on $\bar{\Omega}$. i.e. 
\be 
\Lambda:=\max_{x\in\bar{\Omega}}\{|f(x)|, |Df(x)|, |D^2f(x)|\}
\ene  Let $\Sigma$ be the graph of $u(x)$ with its upward normal vector $\vec{v}=\F{\P_r-Du}{\omega}$ in $M_0$. By \eqref{eq:mc:def}, the mean curvature of $\Sigma$ with respect to $\vec{v}$, $H_\Sigma$, in $M_0$ satisfies that $H_\Sigma+\la Df, \vec{v}\ra=0$. Let $\Theta=\la \vec{v},\P_r\ra=\F{1}{\omega}$ which is a $C^{2,\beta}$ function. By Lemma \ref{lm:angle}, we have 
\be 
\begin{split}
	\Delta \Theta &+(|A|^2+\bar{R}ic(\vec{v},\vec{v}) -Hess(f)(\vec{v},\vec{v}\ra)\Theta \\
	&- \la \nb \Theta, Df\ra =0      
\end{split}
\ene 
Thus $\omega$ satisfies that  
\be \label{eq:estimate:omega}
\F{\Delta \omega}{\omega}\geq  2\F{|\nb \omega|^2}{\omega^2}-c_1+\la \F{\nb \omega}{\omega}, Df \ra 
\ene where $c_1\geq 0$ is a constant only depending on the Ricci curvature and $\Lambda$.  
\\
\indent   Let $d(x, x_0)$ be the distance function between $x_0$ and $x$. Now we define 
\be \label{def:q}
q(x):= 1+\F{u(x)}{2c_0}-\F{3}{2\rho^2}d^2(x,x_0)
\ene 
Define $B:=\{x\in N: q(x)>0\}$. Thus $B_{\F{\rho}{2}}(x_0)\subset B\subset B_{\rho}(x_0)$.  Set $\eta(x) :=e^{K q(x)}-1$ where $K$ is a positive constant determined later. Thus the maximum of $\eta \omega$ is obtained at a point in $B$, for example $x_1\in B_{\rho}(x_0)$. At this point, $\eta\nb \omega+\omega \nb \eta=0$ and 
\begin{align}
0&\geq \F{1}{\omega}\Delta(\eta\omega)=\F{\eta}{\omega}\Delta \omega+2\F{\la\nb \eta, \nb\omega\ra }{\omega}+\Delta \eta\notag\\
&\geq -c_1\eta+\eta \la \F{ \nb\omega}{\omega},Df\ra +\Delta\eta \text{\,by\,}\eqref{eq:estimate:omega}\notag\\
&\geq c_1-c_1 e^{Kq(x)}+\la \nb \eta, Df\ra +\Delta \eta\label{eq:step:es:one}
\end{align}
Observe that 
\be \label{eq:step:es:two}
\la \nb \eta, Df \ra \geq e^{Kq(x)}(-\F{K^2|\nb q(x)|^2}{2}-4\Lambda^2 )
\ene 
Since $\Delta \eta=K\Delta q(x)+K^2|\nb q(x)|^2$, from \eqref{eq:step:es:one} and \eqref{eq:step:es:two} we obtain 
\begin{align}\label{deq:sd}
0&\geq e^{Kq(x)}(K\Delta q(x)+\F{1}{2}K^2 |\nb q(x)|^2-4\Lambda^2-c_1)
\end{align}
Note that for any $C^2$ function $h(x)$, $\Delta h=Hess(h)(e_i,e_i)-H_{\Sigma}\la Dh,\vec{v}\ra)$ where $\{e_i\}_{i=1}^n$ is the orthonormal frame on $S$. Thus 
\begin{gather}\notag
\Delta u= \F{1}{\omega}\la Df, -\F{Du}{\omega}\ra,\\
 \Delta d^2(x,x_0)= Hess(d^2)(e_i, e_i)+\la Df, -\F{Du}{\omega}\ra\la D d^2(x,x_0), -\F{Du}{\omega}\ra \notag
\end{gather}
 As a result, 
\be\label{dedt}
\Delta q(x)\geq  c_2
\ene 
where $c_2$ is  a  constant depending on $\Lambda$, the Ricci curvature on $B_{\rho}(x_0)$. Let $D$ be the covariant derivative on $M$. Note that 
\be \label{dedt2}
|\nb q(x)|^2\geq \F{1}{4c_0}|Du|^2-c_3
\ene 
 Here $c_3$ is a fixed constant depending on $\rho$. Now combining \eqref{deq:sd} with \eqref{dedt}, \eqref{dedt2} together we obtain 
\be 
0\geq \F{c_3}{K}+\F{1}{4 c_0}|Du|^2-c_2-\F{1}{K^2}(4\Lambda^2-c_1)
\ene 
When taking $K$ sufficiently large, we obtain at $x_1$, $|Du|^2\leq \F{1}{2}c_3$. Note that $K$ only depends on $\Lambda$, $c_1$ and $c_2$, the Ricci curvature of $M$ on $B_{\rho}(x_0)$. Thus 
\be 
(e^{Kq(x)}-1)\omega \leq e^{\F{3}{2}K}(1+\F{c_3}{2})
\ene 
on $B_{\rho}(x_0)$. Note that for any $x\in B_{\F{\rho}{2}}(x_0)$,  by \eqref{def:q} $q(x) \geq \F{1}{8}$.  Thus $\max_{B_{\F{\rho}{2}}(x_0)}|Du|\leq C$. Here $C$ is  a constant only depending on $c_0$,$\Lambda$ and $\rho$. \\
\indent The proof is complete. 

\ep 
\bibliographystyle{abbrv}
\bibliography{Ref}
\end{document}